\documentclass[12pt]{amsart}
\usepackage{amssymb}
\usepackage{graphicx}
\newtheorem{Def}{Definition}

\newtheorem{Thm}{Theorem}
\newtheorem{Cor}{Corollary}

\newenvironment{Pf}{ Proof.}{\(\square\)}

\title[On the divergence representation ...]{On the divergence representation of the Gauss curvature of Riemannian surfaces and its applications}
\author{Cs. Vincze, M. Ol\'{a}h and L. M. Alabdulsada}
\address{Inst. of Math., Univ. of Debrecen \\
H-4002 Debrecen, P.O.Box 400 \\
Hungary}
\email{csvincze@science.unideb.hu}
\email{olma4000@gmail.com, layth.muhsin@science.unideb.hu}
\begin{document}
\begin{abstract}In the paper we consider Riemannian surfaces admitting a global expression of the Gauss curvature as the divergence of a vector field. It is equivalent to the existence of a metric linear connection of zero curvature. Such a linear connection $\nabla$ plays an important role in the differential geometry of non-Riemannian surfaces in the sense that the Riemannian quadratic forms can be changed into Minkowski functionals in the tangent planes such that the Minkowskian length of the tangent vectors is invariant under the parallel translation with respect to $\nabla$ (compatibility condition). A smoothly varying family of Minkowski functionals in the tangent planes is called a Finslerian metric function under some regularity conditions. Especially, the existence of a compatible linear connection provides the Finsler surface to be a so-called generalized Berwald surface. It is an alternative of the Riemannian geometry for $\nabla$. Using some general observations and topological obstructions we concentrate on explicit examples. In some representative cases (Euclidean plane, hyperbolic plane etc.) we solve the differential equation of the parallel vector fields to construct a smoothly varying family of Minkowski functionals in the tangent planes such that the Minkowskian length of the tangent vectors is invariant under the parallel translation.
\end{abstract}
\maketitle

\footnotetext[1]{Cs. Vincze is supported by the EFOP-3.6.1-16-2016-00022 project. The project is co-financed by the European Union and the European Social Fund. M. Ol\'{a}h is supported by the \'{U}NKP-18-2 New National Excellence Program of the Ministry of Human Capacities, Hungary.}

\section*{Introduction}

The concept of generalized Berwald manifolds goes back to V. Wagner \cite{Wag1}. They are Finsler manifolds admitting linear connections such that the parallel transports preserve the Finslerian length of tangent vectors (compatibility condition). By the fundamental result of the theory \cite{V5} such a linear connection must be metrical with respect to the averaged Riemannian metric given by integration of the Riemann-Finsler metric on the indicatrix hypersurfaces. Therefore the linear connection is uniquely determined by its torsion tensor. The torsion tensor has a special decomposition in 2D because of 
\begin{equation} 
\label{tor}
T(X,Y)=\left(X^1Y^2-X^2Y^1\right)\left(T_{12}^1 \frac{\partial}{\partial u^1}+T_{12}^2 \frac{\partial}{\partial u^2}\right)=\rho(X)Y-\rho(Y)X,
\end{equation}
where $\rho_1=T_{12}^2$ and $\rho_2=-T_{12}^1=T_{21}^1$. In higher dimensional spaces such a linear connection is called semi-symmetric. Using some previous results \cite{V6}, \cite{V9}, \cite{V10} and \cite{V11}, the torsion tensor of a semi-symmetric compatible linear connection can be expressed in terms of metrics and differential forms given by averaging independently of the dimension of the space, but the compatible linear connection must be of zero curvature in 2D unless the manifold is Riemannian \cite{V12}. Therefore we can conclude some topological obstructions because the existence of a metric linear connection of zero curvature is equivalent to the divergence representation of the Gauss curvature of the Riemannian surface. We prove, for example, that any compact generalized Berwald surface without boundary must have zero Euler characteristic. Therefore the Euclidean sphere does not carry such a geometric structure. An important consequence is that the local conformal flatness is taking to fail in the non-Riemannian differential geometry of surfaces. In some further representative cases (Euclidean plane, hyperbolic plane etc.) we solve the differential equation of the parallel vector fields to present explicit examples of non-Riemannian two-dimensional generalized Berwald manifolds as well.

\section{Notations and terminology}

Let $M$ be a connected differentiable manifold with local coordinates $u^1, \ldots, u^n.$ The induced coordinate system of the tangent manifold $TM$ consists of the functions $x^1, \ldots, x^n$ and $y^1, \ldots, y^n$. For any $v\in T_pM$, $x^i(v)=u^i\circ \pi (v)$ and $y^i(v)=v(u^i)$, where $i=1, \ldots, n$ and $\pi\colon TM\to M$ is the canonical projection. 

A \emph{Finsler metric} is a continuous function $F\colon TM\to \mathbb{R}$ satisfying the following conditions:
$\displaystyle{F}$ is smooth on the complement of the zero section (regularity), $\displaystyle{F(tv)=tF(v)}$ for all $\displaystyle{t> 0}$ (positive homogenity) and the Hessian $\displaystyle{g_{ij}=\frac{\partial^2 E}{\partial y^i \partial y^j}}$, where $E=\frac{1}{2}F^2$, is positive definite at all nonzero elements $\displaystyle{v\in T_pM}$ (strong convexity). The so-called \emph{Riemann-Finsler metric} $g$ is constituted by the components $g_{ij}$. It is defined on the complement of the zero section. The Riemann-Finsler metric makes each tangent space (except at the origin) a Riemannian manifold with standard canonical objects such as the {\emph {volume form}} $\displaystyle{d\mu=\sqrt{\det g_{ij}}\ dy^1\wedge \ldots \wedge dy^n}$,
the \emph {Liouville vector field} $\displaystyle{C:=y^1\partial /\partial y^1 +\ldots +y^n\partial / \partial y^n}$  and the {\emph {induced volume form}}
$$\mu=\sqrt{\det g_{ij}}\ \sum_{i=1}^n (-1)^{i-1} \frac{y^i}{F} dy^1\wedge\ldots\wedge dy^{i-1}\wedge dy^{i+1}\ldots \wedge dy^n$$
 on the indicatrix hypersurface $\displaystyle{\partial K_p:=F^{-1}(1)\cap T_pM\ \  (p\in M)}$. 

A linear connection $\nabla$ on the base manifold $M$ is called {\emph{compatible}} to the Finslerian metric if the parallel transports with respect to $\nabla$ preserve the Finslerian length of tangent vectors. Finsler manifolds admitting compatible linear connections are called \emph{generalized Berwald manifolds}. In case of classical Berwald manifolds the compatible linear connection is trosion-free.

\begin{Thm}
\label{heritage} \emph{\cite{V5}} If a linear connection on the base manifold is compatible with the Finslerian metric function then it must be metrical with respect to the averaged Riemannian metric
\begin{equation}
\label{averagemetric1}
\gamma_p (v,w):=\int_{\partial K_p} g(v, w)\, \mu=v^i w^j \int_{\partial K_p} g_{ij}\, \mu \ \ (v, w\in T_p M, p\in U).
\end{equation}
\end{Thm}

\section{The divergence representation of the Gauss curvature} Let $\nabla$ be a linear connection on the base manifold $M$ of dimension $2$ and suppose that the parallel transports preserve the Finslerian length of tangent vectors (compatibility condition). Theorem \ref{heritage} implies that it is uniquely determined by its torsion tensor of the form
\begin{equation}
\label{torsionform}
T(X,Y)=\rho(X)Y-\rho(Y)X;
\end{equation}
see Formula (\ref{tor}). The idea of the comparison of $\nabla$ with the L\'{e}vi-Civita connection $\nabla^*$ associated with the averaged Riemannian metric (\ref{averagemetric1}) was used to solve the problem of the intrinsic characterization of the semi-symmetric compatible linear connections for both low and higher dimensional spaces. The solution is the expression of the $1$ - form $\rho$ in terms of the canonical data (metrics and differential forms given by averaging) of the Finsler manifold. For the details see  \cite{V6}, \cite{V9}, \cite{V10} and \cite{V11}. Let a point $p\in M$ be given and consider the orthogonal group with respect to the averaged Riemannian metric. It is clear that the subgroup $G\subset O(2)$ of the orthogonal transformations leaving the Finslerian indicatrix invariant is finite unless the Finsler surface reduces to a Riemannian one; see \cite{V12}. If $\nabla$ is a linear connection on the base manifold such that the parallel transports preserve the Finslerian length of tangent vectors (compatibility condition) then, by Theorem \ref{heritage}, $\textrm{Hol}_p \nabla \subset G$ is also finite for any $p\in M$ and the curvature tensor of $\nabla$ is zero. In what follows we are going to compute the relation between the curvatures of $\nabla$ and $\nabla^*$. Taking vector fields with pairwise vanishing Lie brackets on the neighbourhood $U$ of the base manifold, the Christoffel process implies that
\begin{equation}
\label{Cproc}
\gamma(\nabla^*_X Y,Z)=\gamma(\nabla_X Y, Z)+\frac{1}{2}\left(\gamma(X, T(Y,Z))+\gamma(Y, T(X,Z))-\gamma(Z, T(X,Y))\right),
\end{equation}
where $\nabla^*$ denotes the L\'{e}vi-Civita connection. If the torsion is of the form (\ref{torsionform}) then we have that
\begin{eqnarray}
\nabla^*_X Y =  \nabla_X Y +  \rho (Y) X  -  \gamma ( X, Y ) \rho^{\sharp} \ \ \Rightarrow \ \ \nabla_X Y =  \nabla^*_X Y -  \rho (Y) X  +  \gamma ( X, Y ) \rho^{\sharp}, \label{02}
\end{eqnarray}
where $\rho^{\sharp}$ is the dual vector field of $\rho$ defined by $\displaystyle{\gamma(\rho^{\sharp}, X)=\rho(X)}$. 
Consider the curvature tensor
\begin{eqnarray}
R(X,Y)Z= \nabla_{X} \nabla_{Y} Z - \nabla_{Y} \nabla_{X} Z \label{03}
\end{eqnarray}
of $\nabla$. By substituting (\ref{02}) into (\ref{03})
\begin{eqnarray*}
R(X,Y)Z= \nabla_{X} \left(\nabla^*_Y Z -  \rho (Z) Y  + \rho^{\sharp} \gamma ( Z, Y ) \right)-\nabla_{Y} \left(\nabla^*_X Z -  \rho (Z) X  + \rho^{\sharp} \gamma ( Z, X ) \right).
\end{eqnarray*}
Some further direct computations show that 
$$R(X,Y)Z= R^*(X,Y)Z+$$
$$\bigg(\gamma(X,Z) \|\rho^{\sharp}\|^2-\rho(X)\rho(Z)-\left(\nabla_X^* \rho\right)(Z)\bigg)Y+\gamma(Y,Z)\nabla^*_X \rho^{\sharp}+\gamma(Y,Z)\rho(X)\rho^{\sharp}+$$
$$\bigg(\left(\nabla_Y^* \rho \right)(Z)+\rho(Y)\rho(Z)-\gamma(Y,Z) \|\rho^{\sharp}\|^2\bigg) X-\gamma(X,Z)\nabla^*_Y \rho^{\sharp}-\gamma(X,Z)\rho(Y)\rho^{\sharp}.$$
Since the holonomy group of $\nabla$ must be finite in case of a non-Riemannian generalized Berwald surface, we have that $R(X,Y)Z=0$. Taking an orthonormal frame $\gamma(X,Y)=0$, $\gamma(X,X)=\gamma(Y,Y)=1$ at the point of $p\in M$ it follows that 
$$0=\gamma (R^*(X,Y)Y,X)+\rho^2(X)+\rho^2(Y)-\|\rho^{\sharp}\|^2+\gamma \left(\nabla^*_X \rho^{\sharp},X \right)+\left(\nabla^*_Y \rho \right)(Y),$$
where $\displaystyle{\rho^2(X)+\rho^2(Y)-\|\rho^{\sharp}\|^2=0}$ and
$$\left(\nabla^*_Y \rho \right)(Y)=Y\rho(Y)-\rho(\nabla^*_Y Y)=Y\gamma(\rho^{\sharp}, Y)-\rho(\nabla^*_Y Y)=\gamma(\nabla^*_Y \rho^{\sharp}, Y).$$
Therefore
\begin{equation}
\label{cur}
0=\kappa^*(p)+\textrm{div}^* \rho^{\sharp}(p) \ \ \Rightarrow \ \ \kappa^*=-\textrm{div}^* \rho^{\sharp},
\end{equation}
where $\kappa^*$ is the Gauss curvature of the manifold with respect to the averaged Riemannian metric and $\displaystyle{\textrm{div}^* \rho^{\sharp}:=\gamma \left(\nabla^*_X \rho^{\sharp},X \right)+\gamma \left(\nabla^*_Y \rho^{\sharp},Y \right)}$
is the divergence operator. Equation (\ref{cur}) is called \emph{the divergence representation of the Gauss curvature}.

\begin{Cor} A Riemannian surface admits a metric linear connection of zero curvature if and only if its Gauss curvature can be represented as a divergence of a vector field.
\end{Cor}

\begin{Cor} If $M$ is a compact generalized Berwald surface without boundary then it must have zero Euler characteristic.
\end{Cor}
\begin{Pf}
Taking the integral of the divergence representation (\ref{cur}) we have the zero Euler characteristic due to the Gauss-Bonnet theorem and the divergence theorem.
\end{Pf}

\begin{Cor}
\label{forbidden} A two-dimensional Euclidean sphere could not carry Finslerian structures admitting compatible linear connections. 
\end{Cor}

\subsection{Exact and closed Wagner manifolds} It is well-known that any Riemannian surface is locally conformally flat by the (local) solution of the second order elliptic partial differential equation $\displaystyle{\Delta^* f=\kappa^*}$. Its Finslerian analogue is that \emph{any non-Riemannian Finsler surface is locally conformal to a locally Minkowski manifold of dimension $2$}; a locally Minkowski manifold is a Berwald manifold (torsion-free case, i.e. the compatible linear connection is $\nabla^*$) such that $R^*=0$. The solution of the so-called Matsumoto's problem \cite{V6}, see also \cite{V7}, proves that the statement is false in the non-Riemannian Finsler geometry. 
\begin{itemize}
\item[Step 1] By Hashiguchi and Ichyjio's classical theorem \cite{HY}, see also \cite{V1} and \cite{V2}, a Finsler manifold is a conformally Berwald manifold if and only if there exists a semi-symmetric compatible linear connection  with an exact $1$-form $\rho$ in the torsion (\ref{torsionform}). Especially, it is the exterior derivative of the logarithmic scale function $\alpha$ between the (conformally related, see \cite{H2}) Finslerian fundamental functions $\displaystyle{\tilde{F}=e^{\alpha\circ \pi}F}$ up to a minus sign. 
\end{itemize}
\begin{Def} Generalized Berwald manifolds admitting compatible semi-symmetric linear connections with an exact $1$-form $\rho$ in the torsion \eqref{torsionform} are called exact Wagner manifolds. Generalized Berwald manifolds admitting compatible semi-symmetric linear connections with a closed $1$-form $\rho$ in the torsion \eqref{torsionform} are called closed Wagner manifolds.
\end{Def}

\begin{itemize}
\item[Step 2] The generalization of Hashiguchi and Ichyjio's classical theorem for closed Wagner manifolds is the statement that \emph{a Finsler manifold is a locally conformally Berwald manifold if and only if it is a closed Wagner manifold}. It is clear from the global version of the theorem that any point of a closed Wagner manifold has a neighbourhood over which it is conformally equivalent to a Berwald manifold, i.e. any closed Wagner manifold is a locally conformally Berwald manifold. What about the converse? Suppose that we have a locally conformally Berwald manifold. The exterior derivatives of the local scale functions constitute  a globally well-defined closed $1$-form for the torsion (\ref{torsionform}) of a compatible linear connection if and only if they coincide on the intersections of overlapping neighbourhoods. Since the conformal equivalence is transitive it follows that overlapping neighbourhoods carry conformally equivalent Berwald metrics. The problem posed by M. Matsumoto \cite{M2} in 2001 is that \emph{are there non-homothetic and non-Riemannian conformally equivalent Berwald spaces}? It has been completely solved by \cite{V6} in 2005, see also \cite{V7}.
\end{itemize}

\begin{Thm} \emph{\cite{V6}, see also \cite{V7}} The scale function between conformally equivalent Berwald manifolds must be constant unless they are Riemannian. 
\end{Thm}

\begin{Cor} \emph{\cite{V6}, see also \cite{V7}} A Finsler manifold is a locally conformally Berwald manifold if and only if it is a closed Wagner manifold.
\end{Cor}

Using Corollary \ref{forbidden} we have the following result.

\begin{Cor} A two-dimensional Euclidean sphere could not carry non-Riemannian locally conformally Berwald Finslerian structures. Especially, it can not be a locally conformally flat non-Riemmanian Finsler manifold.
\end{Cor}

By the classification of orientable compact surfaces without boundary we can also state that they could not carry Finslerian structures admitting compatible linear connections except the case of genus $1$. In what follows we formulate a result including the generic case of the tori $S^1\times S^1$.

\begin{Thm} The Gauss curvature has a divergence representation for any Riemannian manifold of the form $M=M_1\times M_2$, where $M_i=\mathbb{R}$ or $S^1$, $S^1$ is the Euclidean unit circle and  $i=1, 2$.
\end{Thm}

\begin{Pf}
Let us start the proof with the case of $M=\mathbb{R}^2$. Using formula 
\begin{equation}
\label{divgeneral}
\textrm{div}^* (X) = \frac{1}{\sqrt{\textrm{det} \, \gamma_{ij}}} \left[ \frac{\partial \left( \sqrt{\textrm{det} \, \gamma_{ij}} X^1 \right)}{\partial u^1} + \frac{\partial \left( \sqrt{\textrm{det} \, \gamma_{ij}} X^2 \right)}{\partial u^2} \right]
\end{equation}
it follows that 
$$X^2(u^1,u^2)=$$
$$-\frac{1}{\sqrt{\det \gamma_{ij}(u^1,u^2)}}\left(\int_{0}^{u^2}\kappa^*(u^1, t) \sqrt{\det \gamma_{ij} (u^1,t)}+\frac{\partial \left( \sqrt{\det \gamma_{ij}} X^1 \right)}{\partial u^1}(u^1,t)\, dt+c(u^1)+c_0\right).
$$
If the functions of the right hand side satisfy the periodicity conditions
\begin{equation}
\begin{gathered}
\det \gamma_{ij}(u^1+u^1_0,u^2+u_0^2)=\det \gamma_{ij}(u^1,u^2), \ \kappa^*(u^1+u^1_0,u^2+u_0^2)=\kappa(u^1,u^2),\\
X^1(u^1+u^1_0,u^2+u_0^2)=X^1(u^1,u^2) \ \textrm{and} \ c(u^1+u^1_0)=c(u^1)\\
\end{gathered}
\end{equation}
then $\displaystyle{X^2(u^1+u^1_0,u^2+u_0^2)=X^2(u^1,u^2)}$. The case $u^1_0=u_0^2=0$ belongs to $M=\mathbb{R}^2$. If $u^1_0=0$ and $u_0^2>0$ then we can interpret $X^1$ and $X^2$ as functions on $M:=\mathbb{R}\times S^1$. Finally, the case $u_0^1>0$ and $u_0^2>0$ corresponds to the solution on the tori $S^1\times S^1$. 
\end{Pf}

\section{The Euclidean and the hyperbolic plane}

In what follows we present more explicit solutions of the divergence representation problem of the curvature. On the other hand we use the dual forms of the vector fields to determine the parallel translations with respect to the metric linear connection $\nabla$ of torsion (\ref{torsionform}). As we shall see, the holonomy groups are trivial. Using parallel translations for the extension, it is enough to substitute the Riemannian indicatrix with a more general convex closed curve\footnote{In the forthcoming examples we will use trifocal ellipses as convex closed curves containing the origin in their interiors.} containing the origin in its interior at a single point. Such a smoothly varying family of curves provides an alternative way of measuring the length of tangent vectors, i.e. we have a non-Riemannian generalized Berwald manifold with $\nabla$ as the compatible linear connection. 

\subsection{The case of the Euclidean plane} Consider the Euclidean plane $\mathbb{R}^2$ equipped with the canonical inner product $\delta_{ij}$. The divergence representation of the Gauss curvature means to find a vector field  $\rho^{\sharp}$ with vanishing divergence. Since the curl of the rotated vector field
$$\rho^2 \frac{\partial}{ \partial u^1}-\rho^1 \frac{\partial}{ \partial u^2}$$
is zero, we have a global solution (potential) of equations $\displaystyle{\rho^2=\frac{\partial f}{ \partial u^1}}$ and $\displaystyle{\rho^1=-\frac{\partial f}{ \partial u^2}}$. 

Taking $\rho_i=\delta_{ik}\rho^k=\rho^i$ ($i=1, 2$) let $\displaystyle{\rho=\rho_1du^1+\rho_2du^2}$ be the dual $1$-form of $\rho^{\sharp}$ and consider the metric linear connection $\nabla$ of torsion (\ref{torsionform}). If $X$ is a parallel vector field with respect to $\nabla$ along a curve $\displaystyle{c\colon [0, 1] \to \mathbb{R}^2}$, then, by formula (\ref{02}),
$$\left(X^1\right)'=\left(c^1\right)'\rho_c (X)-\left(\left(c^1\right)'X^1+\left(c^2\right)'X^2\right)\rho^1\circ c,$$
$$\left(X^2\right)'=\left(c^2\right)'\rho_c (X)-\left(\left(c^1\right)'X^1+\left(c^2\right)'X^2\right)\rho^2\circ c.$$
Therefore
$$\left(X^1\right)'=X^2\left(\left(c^1\right)'\rho_2 \circ c-\left(c^2\right)'\rho_1\circ c\right) \ \ \textrm{and}\ \ \left(X^2\right)'=X^1\left(\left(c^2\right)'\rho_1 \circ c-\left(c^1\right)'\rho_2\circ c \right)$$
and the differential equations of the parallel vector fields are $\displaystyle{\left(X^1\right)'=\varphi' X^2}$ and $\displaystyle{\left(X^2\right)'=-\varphi' X^1}$, where $\displaystyle{\varphi=f\circ c}$. If $\nabla$ is metrical then 
$$X(t)=r_0\left(\cos \theta(t) \frac{\partial}{ \partial u^1}\circ c(t)+\sin \theta(t) \frac{\partial}{ \partial u^2}\circ c(t)\right)$$
with constant Euclidean norm $r_0$, where $\theta'(t)=-\varphi'(t)$ because of the parallelism. The general form of a parallel vector field with respect to $\nabla$ along the curve $\displaystyle{c\colon [0, 1] \to \mathbb{R}^2}$ is
$$X(t)=r_0\left(\cos \left(\varphi(t)+\varphi_0\right) \frac{\partial}{ \partial u^1}\circ c(t)-\sin \left(\varphi(t)+\varphi_0\right) \frac{\partial}{ \partial u^2}\circ c(t)\right).$$
It is clear that if $c(0)=c(1)$, then $X(0)=X(1)$, i.e. the holonomy group of $\nabla$ contains only the identity. Taking an arbitrary convex curve around the origin we can extend it by parallel transports with respect to $\nabla$ to the entire plane $\mathbb{R}^2$. Such a collection of indicatrices constitutes a Finslerian metric function $F$ with $\nabla$ as the compatible linear connection.

\subsubsection{An example} If $\displaystyle{\rho=u^2du^1-u^1du^2}$ then $\displaystyle{f(u^1, u^2)=-\frac{1}{2}\left(\left(u^1\right)^2+\left(u^2\right)^2\right)}$ and the parallel vector fields are of the form
$$X(t)=X^1(t)\frac{\partial}{ \partial u^1}\circ c(t)+X^2(t)\frac{\partial}{ \partial u^2}\circ c(t),$$
where
$$X^1(t)=r_0 \cos \left(\frac{1}{2}\left(\left(c^1\right)^2(t)+\left(c^2\right)^2(t)\right)+\varphi_0\right),$$
$$X^2(t)=-r_0\sin \left(\frac{1}{2}\left(\left(c^1\right)^2(t)+\left(c^2\right)^2(t)\right)+\varphi_0\right).$$
Let the trifocal ellipse defined by 
\begin{equation}
\label{3ell}
\sqrt{(u^1+1)^2+(u^2)^2}+\sqrt{(u^1)^2+(u^2)^2}+\sqrt{(u^1-1)^2+(u^2)^2}=4
\end{equation}
be choosen as the indicatrix at the origin. The focal set contains the elements
$$-X_0:=(-1,0), \ {\bf 0}, \ X_0:=(1,0).$$
The parallel translates of the trifocal ellipse (\ref{3ell}) are given by the equations
\begin{equation}
\label{3elltrans}
\sqrt{(u^1+X^1(t))^2+(u^2+X^2(t))^2}+\sqrt{(u^1)^2+(u^2)^2}+\sqrt{(u^1-X^1(t))^2+(u^2-X^2(t))^2}=4,
\end{equation}
where $X$ is a parallel vector field along a curve $c$ satisfying the initial conditions $X^1(0)=1$ and $X^2(0)=0$. The focal set at the parameter $t$ is $\displaystyle{-X(t), \ {\bf 0}, \ X(t)}$. Figure 1 shows the parallel translates of the indicatrix along the radial direction $c(t)=(t,t)$.
\begin{figure}[h]
\includegraphics[scale=0.18]{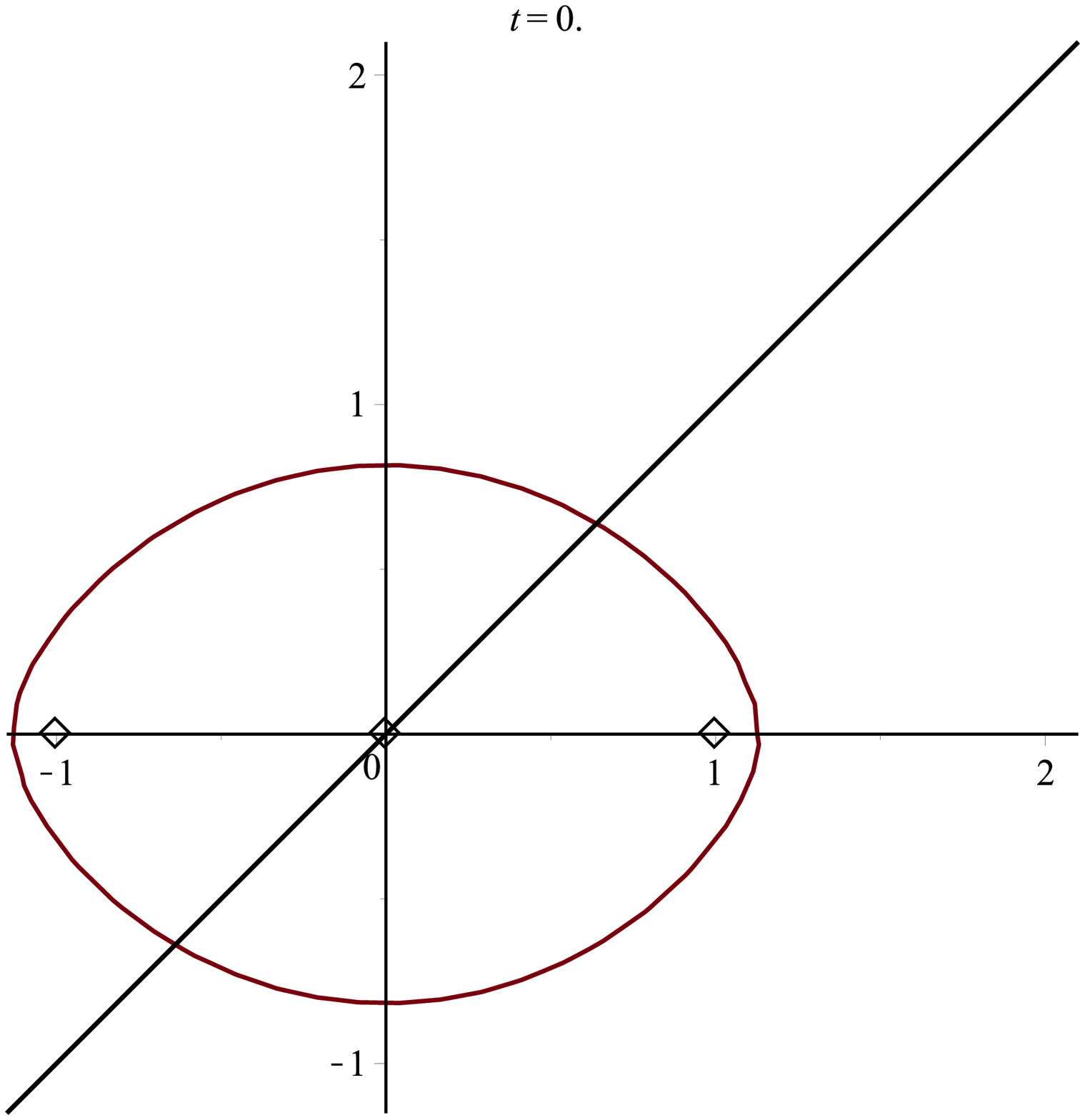}
\includegraphics[scale=0.18]{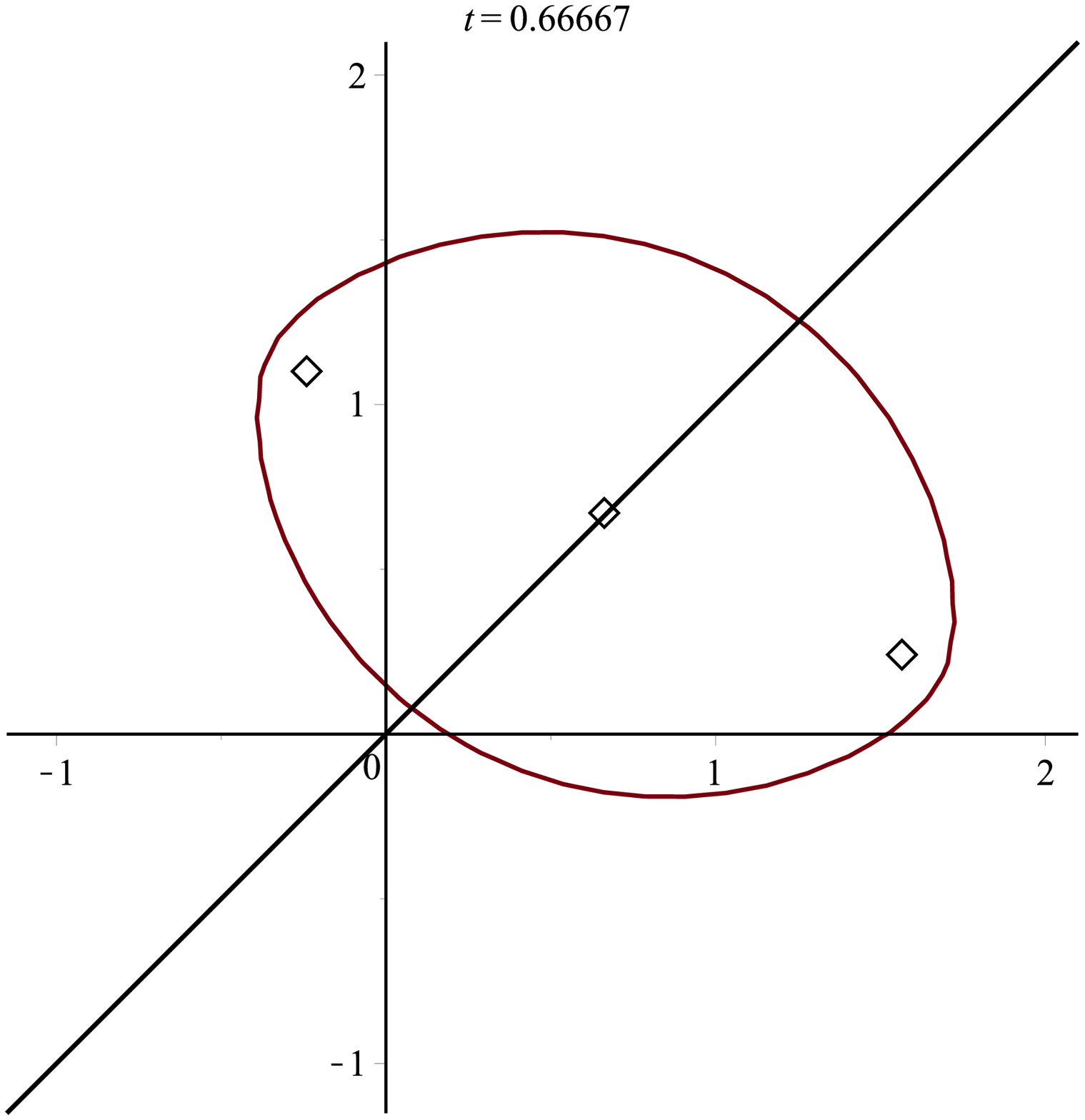}
\includegraphics[scale=0.18]{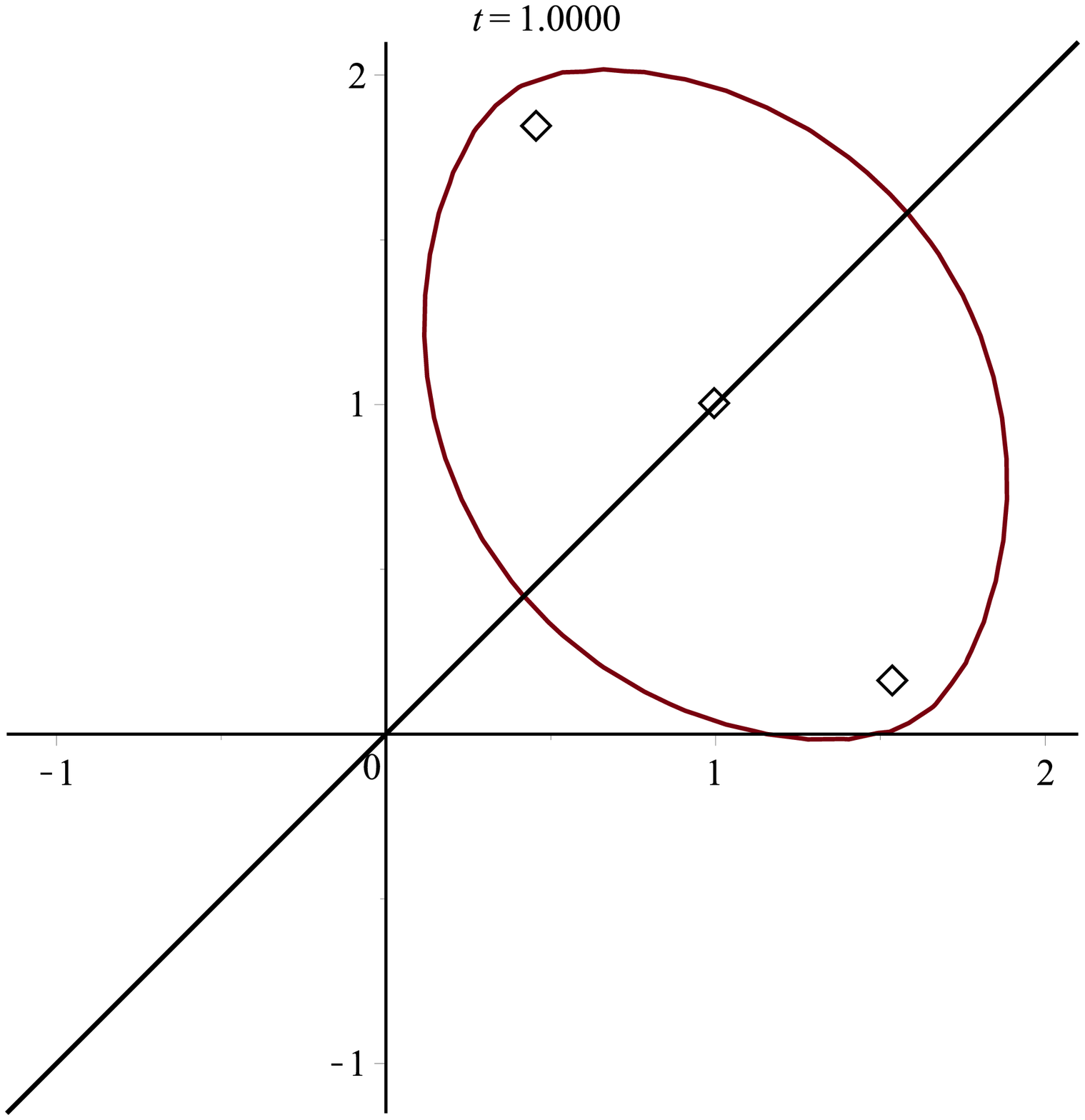}
\caption{Parallel translation along the radial direction (from left to right)}
\end{figure}

In case of the radial direction the focal set of the translated trifocal ellipse at the parameter $t$ is $-X(t)$, ${\bf 0}$,
$$X(t)=\cos (t^2) \frac{\partial}{ \partial u^1}_{(t,t)}-\sin (t^2) \frac{\partial}{ \partial u^2}_{(t,t)}.$$
Figure 2 shows the parallel translates of the indicatrix along the circle $c(t)=(\cos(t), \sin(t)+1)$. The focal set of the translated trifocal ellipse at the parameter $t$ is $-X(t)$, ${\bf 0}$
$$X(t)=\cos \left(1+\sin (t)\right) \frac{\partial}{ \partial u^1}_{(\cos(t),1+\sin(t))}-\sin \left(1+\sin(t)\right) \frac{\partial}{ \partial u^2}_{(\cos(t),1+\sin(t))}.$$
\begin{figure}
\includegraphics[scale=0.18]{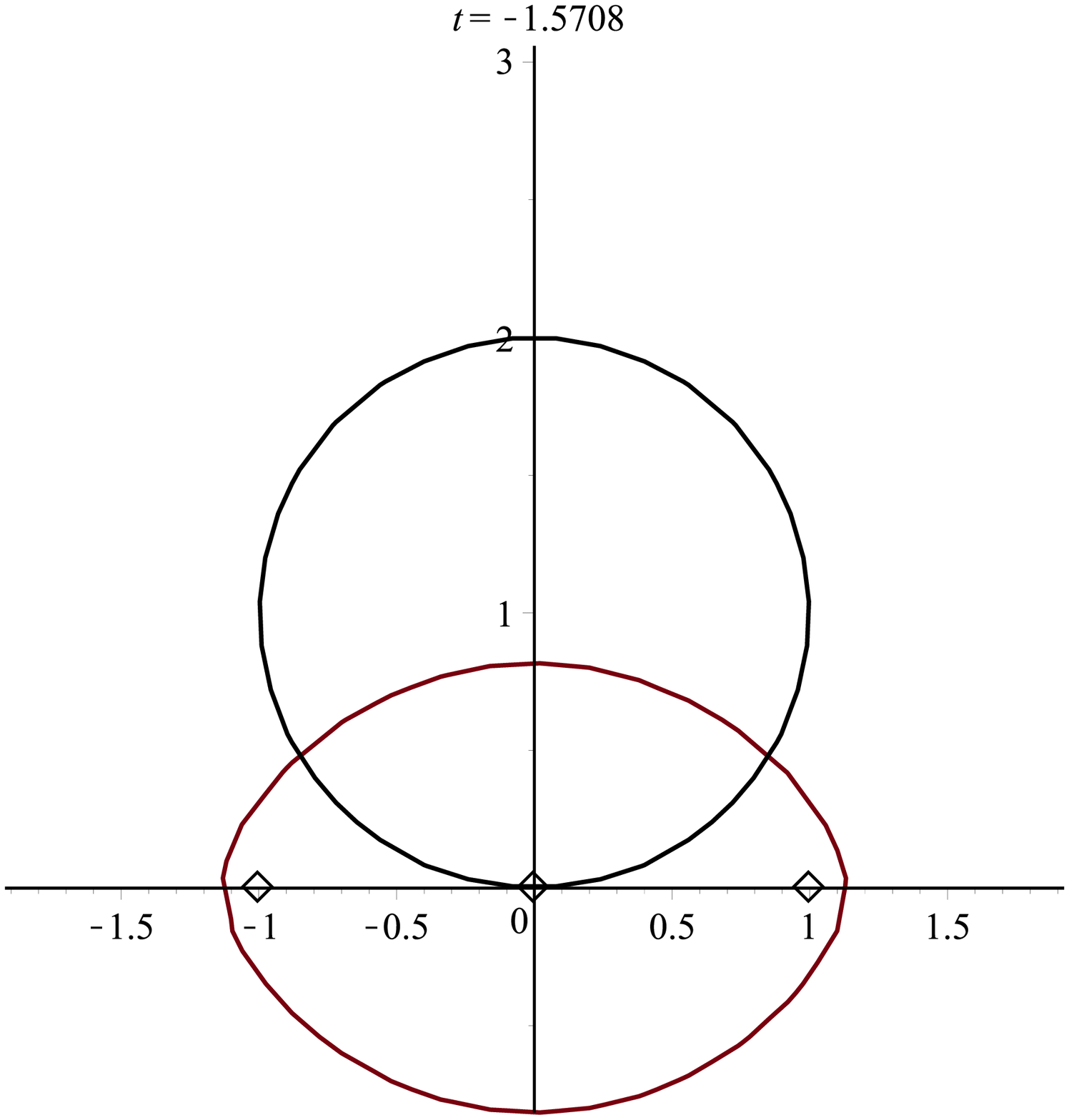}
\includegraphics[scale=0.18]{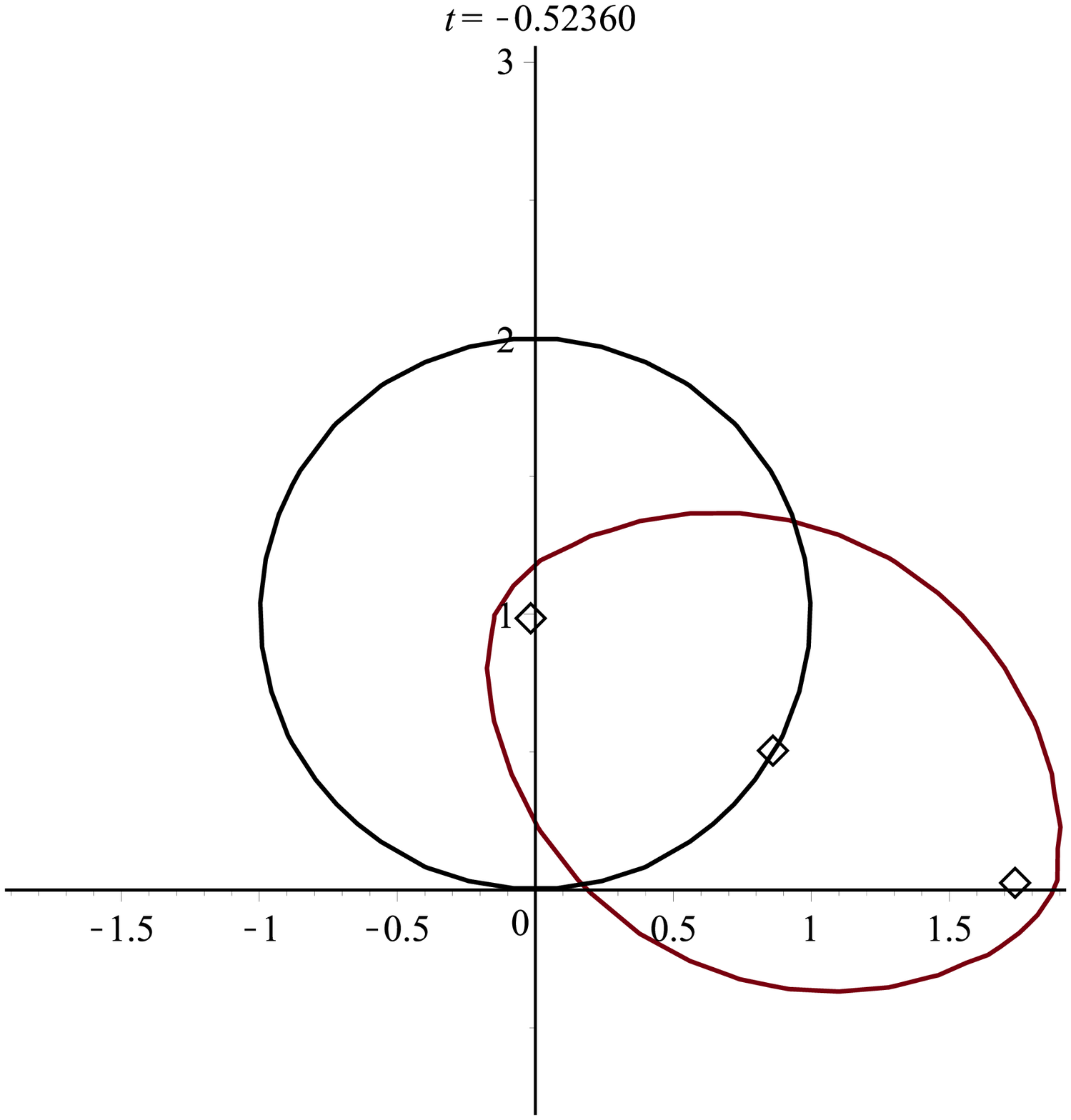}
\includegraphics[scale=0.18]{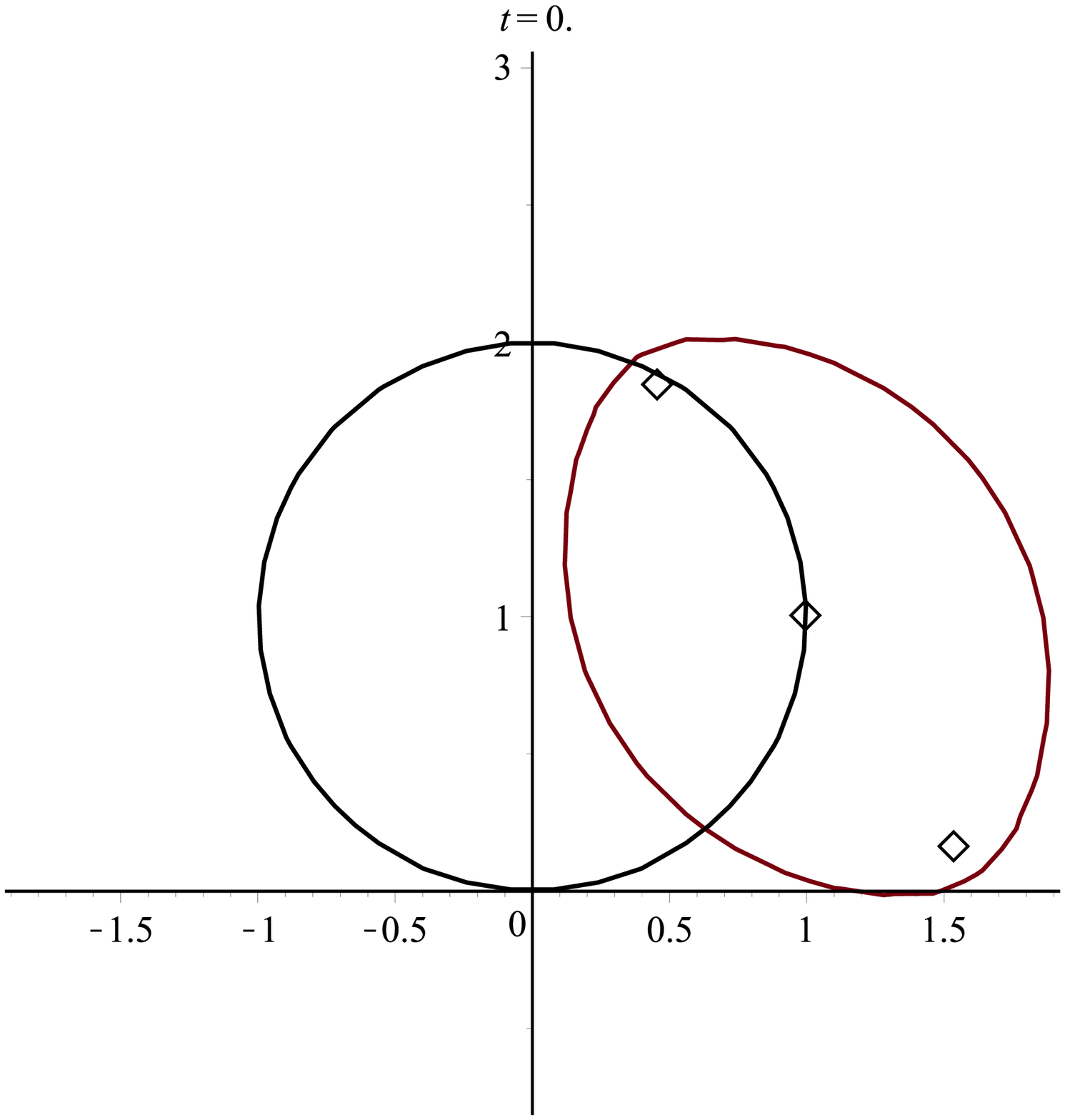}
\end{figure}
\begin{figure}
\includegraphics[scale=0.18]{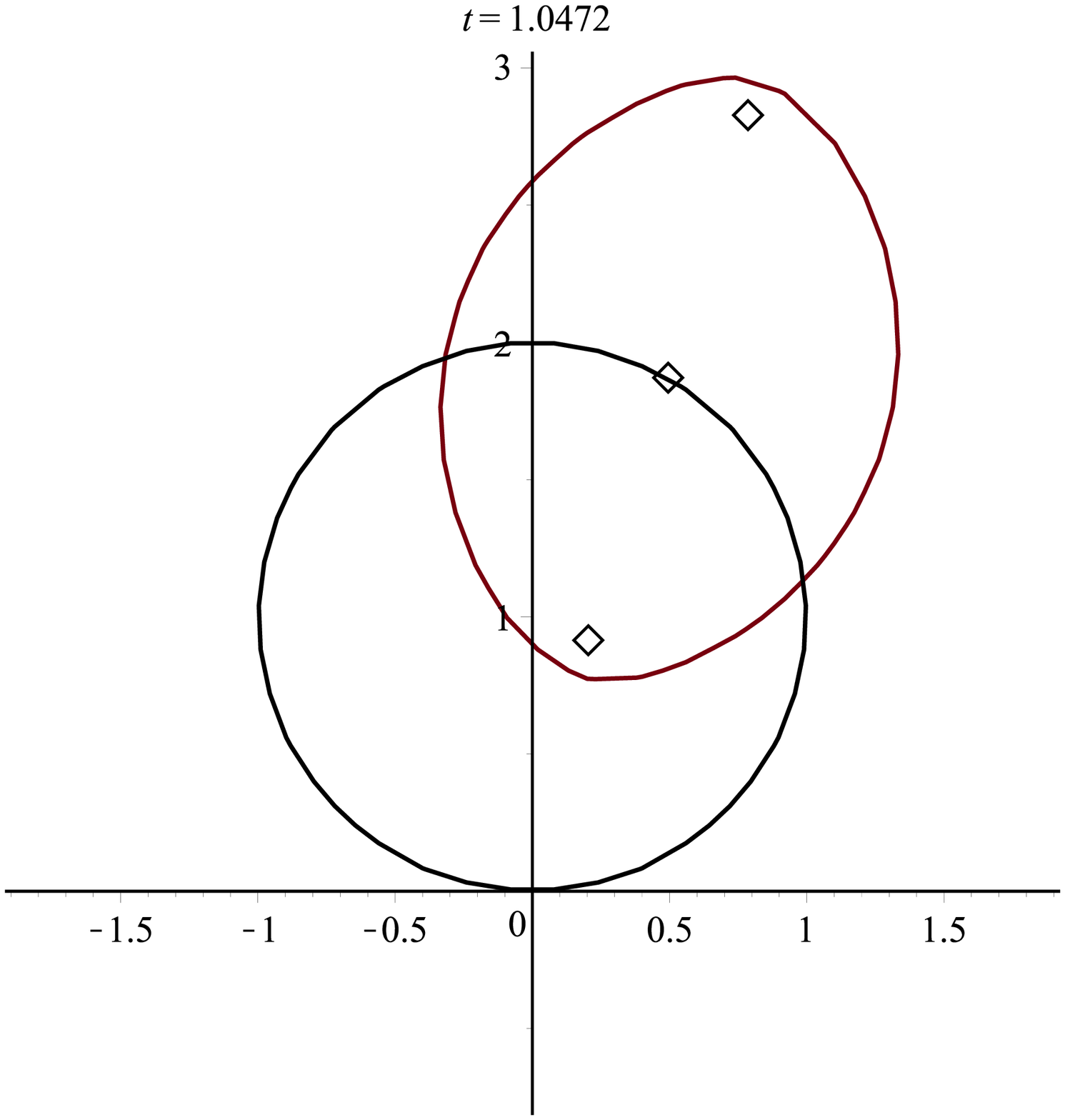}
\includegraphics[scale=0.18]{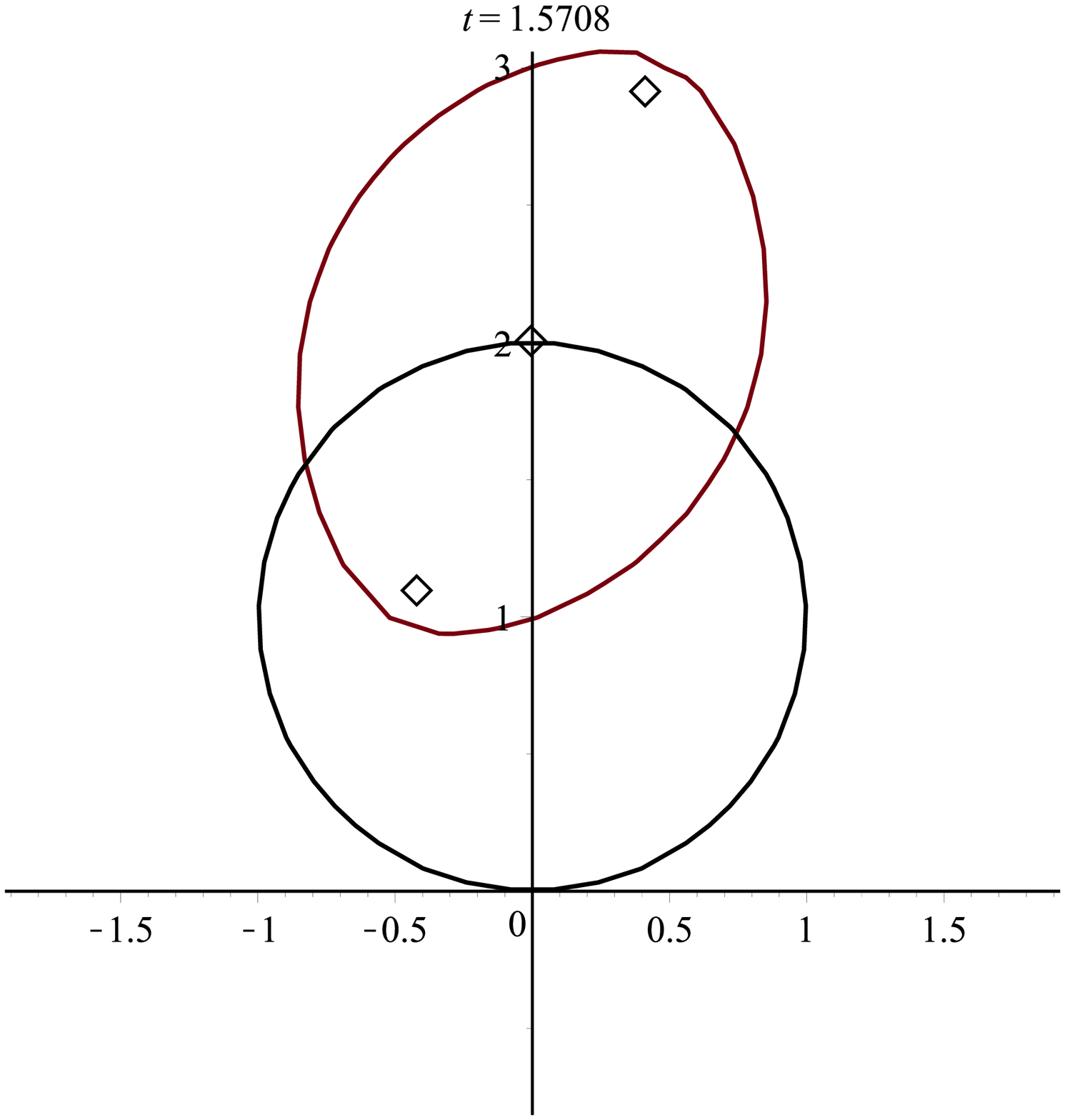}
\includegraphics[scale=0.18]{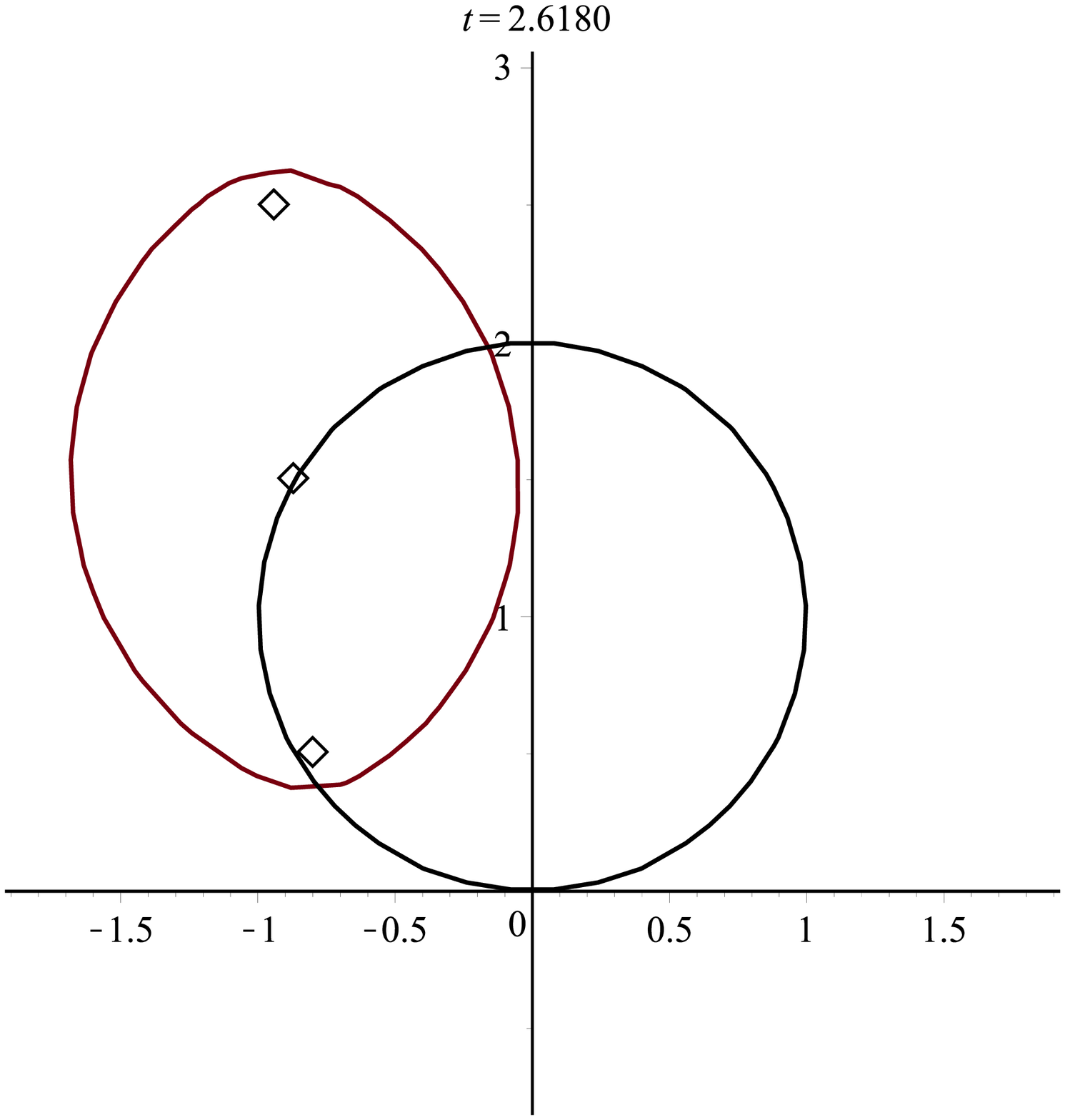}
\end{figure}
\begin{figure}
\includegraphics[scale=0.18]{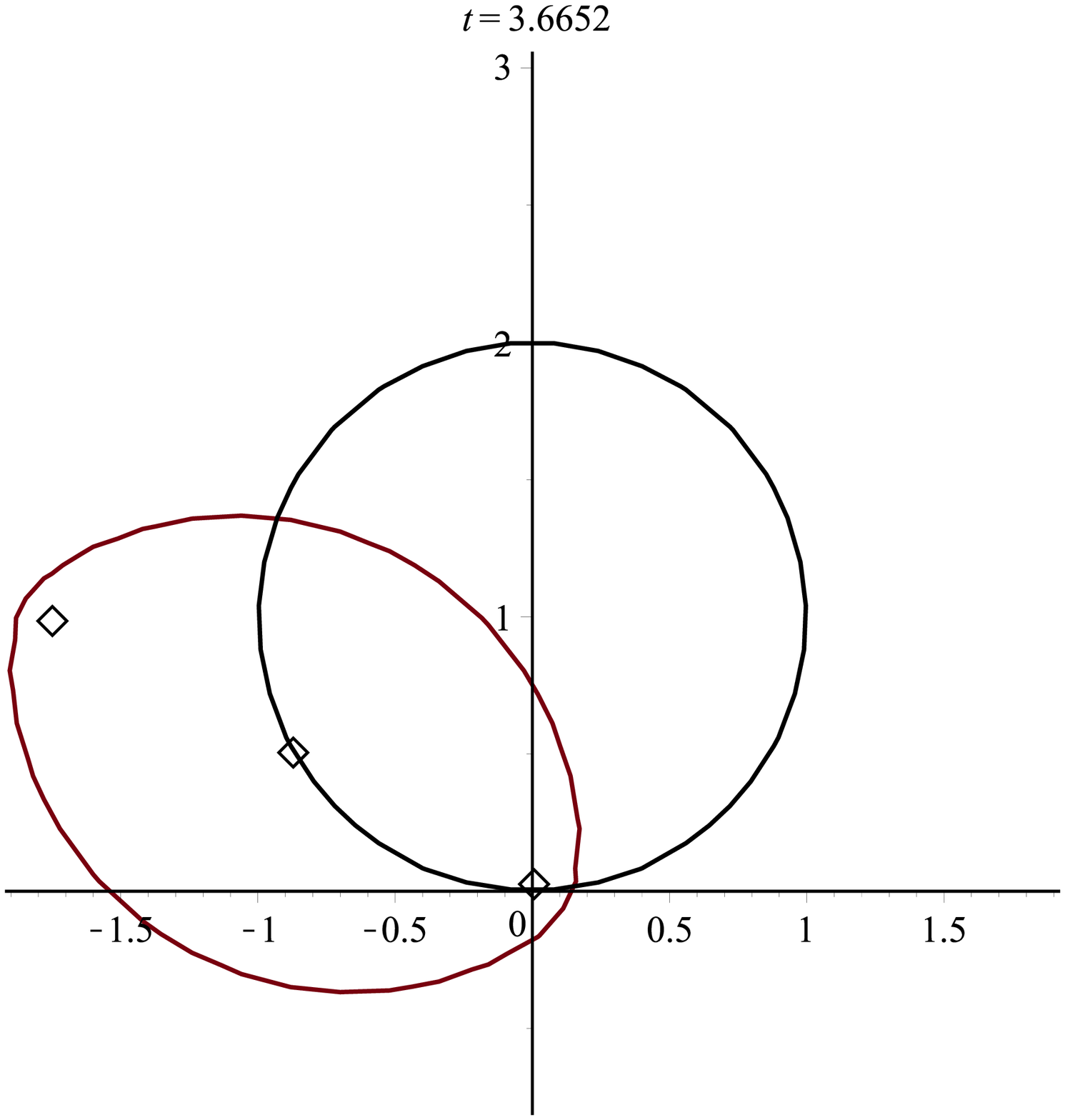}
\includegraphics[scale=0.18]{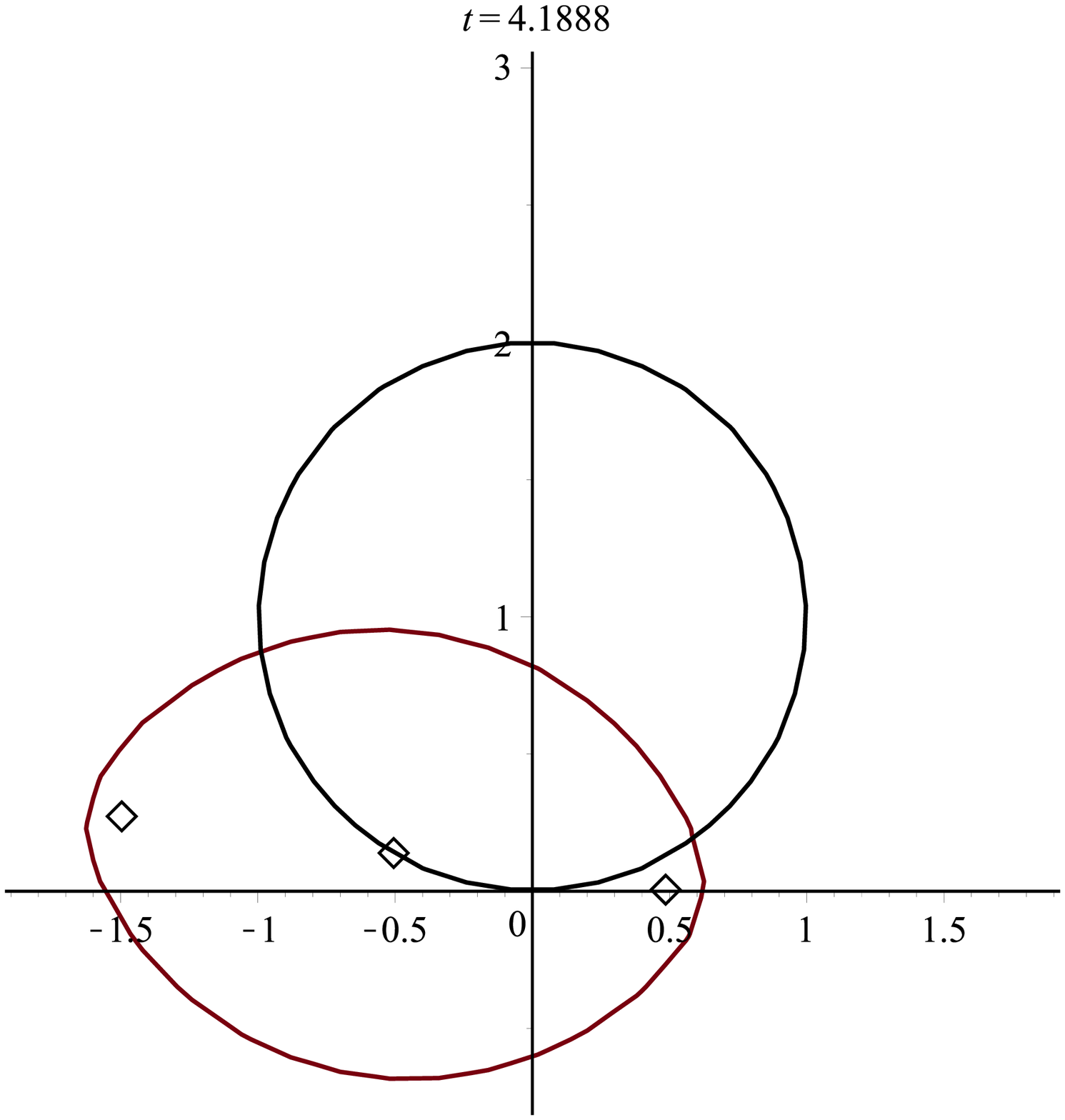}
\includegraphics[scale=0.18]{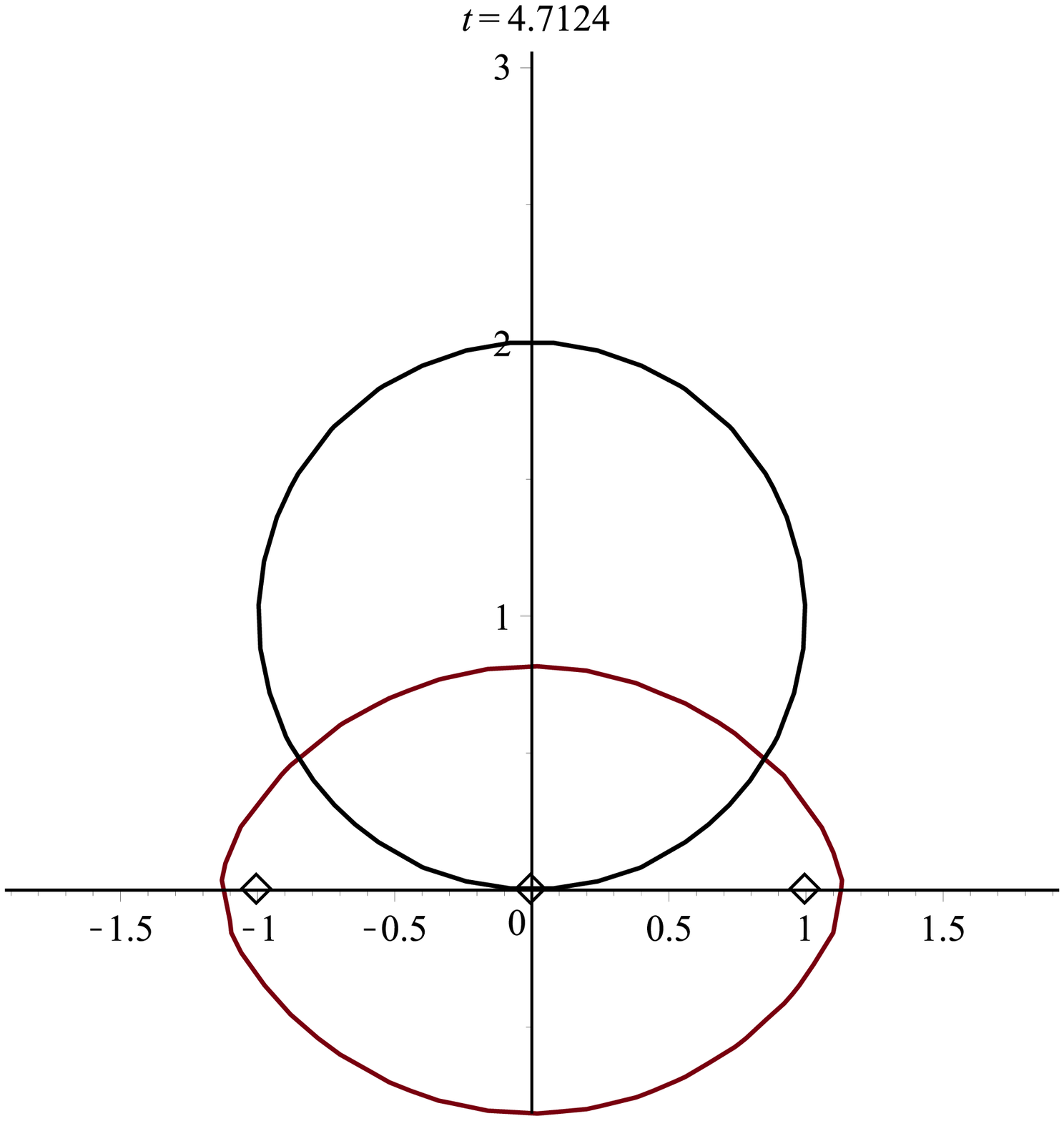}
\caption{Parallel translation along the circle (from left to right)}
\end{figure}

The induced generalized Berwald plane is not conformally flat. To present examples for conformally flat generalized Berwald manifold it is sufficient and necessary to choose a closed (and, consequently exact) $1$-form $\rho$.

\subsection{The case of the hyperbolic plane} Consider the upper half plane 
$$\mathbb{H}_2 := \{ (u^1,u^2) \, | \, u^2 > 0 \}$$
equipped with the metric $\displaystyle{\gamma_{ij} := \frac{1}{(u^2)^2}\delta_{ij}}$. The parameters of the L\'{e}vi-Civita connection $\nabla^*$ are 
\begin{equation}
\label{chrishyp}
\begin{gathered}
\Gamma_{11}^{*1} = 0, \ \Gamma_{12}^{*1} = -\frac{1}{u^2}, \ \Gamma_{21}^{*1} = -\frac{1}{u^2}, \Gamma_{22}^{*1} = 0\\
\Gamma_{11}^{*2} = \frac{1}{u^2}, \ \Gamma_{12}^{*2} = 0 , \ \Gamma_{21}^{*2} = 0, \ \Gamma_{22}^{*2} = -\frac{1}{u^2}. 
\end{gathered}
\end{equation}
The divergence representation (\ref{cur}) of the Gauss curvature means to find a vector field  $\rho^{\sharp}$ such that
$$ -1 = -\textrm{div}^* \rho^{\sharp},$$
where 
$$ \textrm{div}^* \rho^{\sharp} \stackrel{\eqref{divgeneral}}{=} \frac{\partial \rho^1}{\partial u^1} + \frac{\partial \rho^2}{\partial u^2} - \frac{2}{u^2} \rho^2,$$
because of $\displaystyle{\det \, \gamma_{ij} = \frac{1}{(u^2)^4}}$. Therefore 
\begin{equation}\label{div}
 1 = \frac{\partial \rho^1}{\partial u^1} + \frac{\partial \rho^2}{\partial u^2} - \frac{2}{u^2} \rho^2.
\end{equation}
As a straightforward calculation shows, equation (\ref{div}) is equivalent to 
\begin{equation}
0 = \frac{\partial \tilde{\rho}^1}{\partial u^1}+\frac{\partial \tilde{\rho}^2}{\partial u^2},
\end{equation}
where $\displaystyle{\tilde{\rho}^1= \frac{\rho^1}{(u^2)^2}}$ and $\displaystyle{\tilde{\rho}^2= \frac{\rho^2}{(u^2)^2} + \frac{1}{u^2}}$. Using the case of the Euclidean plane,
\begin{equation}
- \tilde{\rho}^2 := - \frac{\rho^2}{(u^2)^2} - \frac{1}{u^2} = \frac{\partial f}{ \partial u^1}, \ \tilde{\rho}^1:= \frac{\rho^1}{(u^2)^2} = \frac{\partial f}{ \partial u^2}.
\end{equation}
The general form of the solution is 
\begin{equation}
\rho^1= (u^2)^2 \frac{\partial f}{ \partial u^2}, \ \rho^2 = - (u^2)^2 \frac{\partial f}{ \partial u^1} - u^2.
\end{equation}
Taking the dual $1$-form $\displaystyle{\rho=\rho_1du^1+\rho_2du^2}$ of $\rho^{\sharp}$, 
$$\rho_i = \gamma_{ik} \rho^k = \frac{1}{(u^2)^2} \rho^i \ (i=1, 2), \ \textrm{i.e.}\ \rho_1= \frac{1}{(u^2)^2} \rho^1 = \frac{\partial f}{ \partial u^2} \ \textrm{and}\ \rho_2= \frac{1}{(u^2)^2} \rho^2 = -\frac{\partial f}{ \partial u^1} - \frac{1}{u^2},$$
consider the metric linear connection $\nabla$ of torsion (\ref{torsionform}). By formula (\ref{02}) and (\ref{chrishyp}) 
\begin{equation}
\label{chris}
\begin{gathered}
\Gamma_{11}^{1}  =  -\rho_1 + \rho_1 = 0, \ \Gamma_{12}^{1} = -\frac{1}{u^2} -\rho_2 = \frac{\partial f}{ \partial u^1}, \ \Gamma_{21}^{1} =  -\frac{1}{u^2}, \ \Gamma_{22}^{1}  =  \rho_1 = \frac{\partial f}{ \partial u^1} \\
\Gamma_{11}^{2} = \frac{1}{u^2} + \rho_2 = -\frac{\partial f}{ \partial u^1}, \ \Gamma_{12}^{2}= 0, \ \Gamma_{21}^{2} = -\rho_1 = -\frac{\partial f}{ \partial u^2}, \ \Gamma_{22}^{2} = -\frac{1}{u^2} -\rho_2 + \rho_2 =  -\frac{1}{u^2}.
\end{gathered}
\end{equation}
If $X$ is a parallel vector field with respect to $\nabla$ along a curve $\displaystyle{c\colon [0, 1] \to \mathbb{H}_2}$, then 
\begin{equation}\label{pareq1}
\begin{gathered}
(X^1)' + (c^i)' X^j  \Gamma_{ij}^1 \circ c  = 0 \\
(X^1)' + (c^1) ' X^2 \frac{\partial f}{ \partial u^1}\circ c -(c^2)' X^1 \frac{1}{u^2}  + (c^2)' X^2 \frac{\partial f}{ \partial u^2}\circ c= 0 \\
(X^1)' + X^2 \left(  (c^1)' \frac{\partial f}{ \partial u^1}\circ c + (c^2)' \frac{\partial f}{ \partial u^2}\circ c \right) -(c^2)' X^1 \frac{1}{c^2} =0 
\end{gathered}
\end{equation}
and, consequently, 
\begin{equation}
\label{pareq11}
(X^1)' + X^2 \varphi' - X^1 \frac{(c^2)'}{c^2} = 0,
\end{equation}
where $\varphi := f \circ c : \mathbb{R} \to \mathbb{R}$. In a similar way, it follows that 
\begin{equation}\label{pareq2}
 (X^2)' - X^1 \varphi' - X^2 \frac{(c^2)'}{c^2} = 0
\end{equation}
To find the solutions we use polar coordinate representations in the tangent planes: taking the global orthonormal frame 
$$ \left( y \frac{\partial}{ \partial x}, y \frac{\partial}{ \partial y} \right) $$
we can write that
\begin{equation}
X(t) = c^2(t) r_0  \left( \cos(\theta(t)) \frac{\partial}{ \partial u^1} \circ c(t) + \sin(\theta(t))\frac{\partial}{ \partial u^2} \circ c(t) \right),
\end{equation}
where $r_0$ is a positive constant becuse $\nabla$ is a metric linear connection. The coordinate functions of $X$ are 
\begin{equation}
X^1= c^2 r_0  \cos (\theta) \ \textrm{and} \ X^2= c^2 r_0 \sin (\theta). \\
\end{equation}
Substituting into \eqref{pareq1} and \eqref{pareq2},
\begin{equation}
\theta' \sin(\theta) = \varphi' \sin(\theta) \ \textrm{and} \ \theta' \cos(\theta)  = \varphi' \cos(\theta).
\end{equation}
We have that 
\begin{equation}
\theta'(t) = \varphi'(t) \ \Rightarrow \ \theta(t) = \varphi(t) + \varphi_0
\end{equation}
and the general form of the parallel vector fields with respect to $\nabla$ is 
\begin{equation}
X(t)=c^2(t) r_0 \left(\cos \left(\varphi(t)+\varphi_0\right) \frac{\partial}{ \partial u^1}\circ c(t) + \sin \left(\varphi(t)+\varphi_0\right) \frac{\partial}{ \partial u^2}\circ c(t) \right).
\end{equation}
It is clear that if $c(0)=c(1)$, then $X(0)=X(1)$, i.e. the holonomy group of $\nabla$ contains only the identity. Taking an arbitrary convex curve at a single point, we can extend it by parallel transports with respect to $\nabla$ to the upper half plane $\mathbb{H}_2$. Such a collection of indicatrices constitutes a Finslerian metric function $F$ with $\nabla$ as the compatible linear connection.

\subsubsection{An example} If $\displaystyle{\rho=\frac{1}{u^2}\left(du^1-du^2 \right)}$, then $\displaystyle{f(u^1, u^2)=\log u^2}$ and the parallel vector fields are of the form
$$X(t)=X^1(t)\frac{\partial}{ \partial u^1}\circ c(t)+X^2(t)\frac{\partial}{ \partial u^2}\circ c(t),$$
where
$$X^1(t)= c^2(t) r_0 \cos \left(\log \left(c^2(t)\right)+\varphi_0\right), \ X^2(t)= c^2(t) r_0 \sin \left(\log \left(c^2(t)\right)+\varphi_0\right).$$
Let the trifocal ellipse defined by 
\begin{equation}
\label{3ellhyp}
\sqrt{(u^1+1)^2+(u^2)^2}+\sqrt{(u^1)^2+(u^2)^2}+\sqrt{(u^1-1)^2+(u^2)^2}=4
\end{equation}
be choosen as the indicatrix at the point $(0,1)\in \mathbb{H}_2$. The focal set contains the elements
$$-X_0:=(-1,0), \ {\bf 0}, \ X_0:=(1,0).$$
The parallel translates of the trifocal ellipse (\ref{3ellhyp}) are given by the equations
\begin{equation}
\label{3elltranshyp}
\sqrt{(u^1+X^1(t))^2+(u^2+X^2(t))^2}+\sqrt{(u^1)^2+(u^2)^2}+\sqrt{(u^1-X^1(t))^2+(u^2-X^2(t))^2}=4,
\end{equation}
where $X$ is a parallel vector field along a curve $c$ satisfying the initial conditions $X^1(0)=1$ and $X^2(0)=0$. The focal set at the parameter $t$ is $\displaystyle{-X(t), \ {\bf 0}, \ X(t)}$. Figure 3 shows the parallel translates of the indicatrix along the straigh line  $c(t)=(t,t+1)$.
\begin{figure}[h]
\includegraphics[scale=0.25]{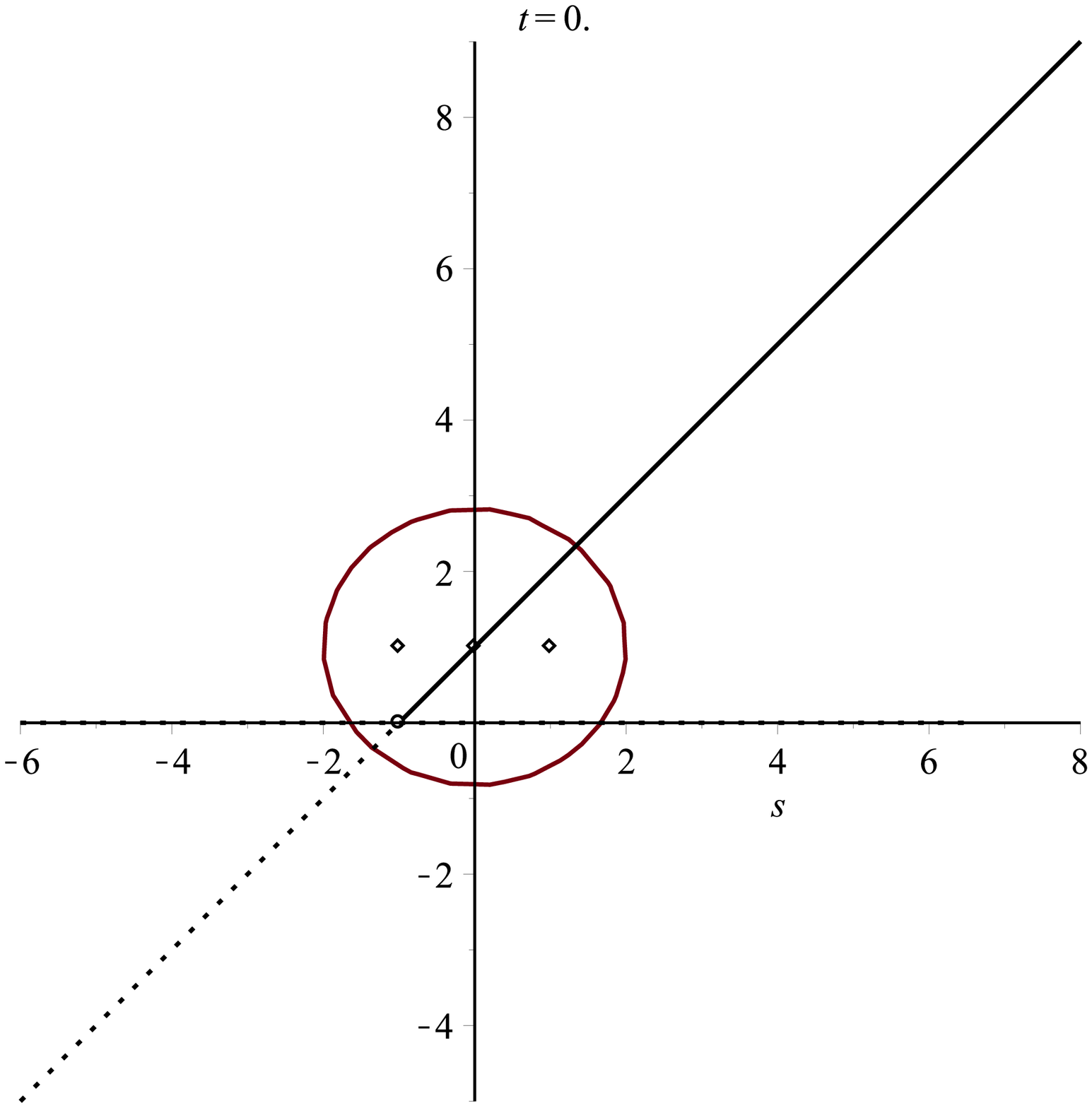}
\includegraphics[scale=0.25]{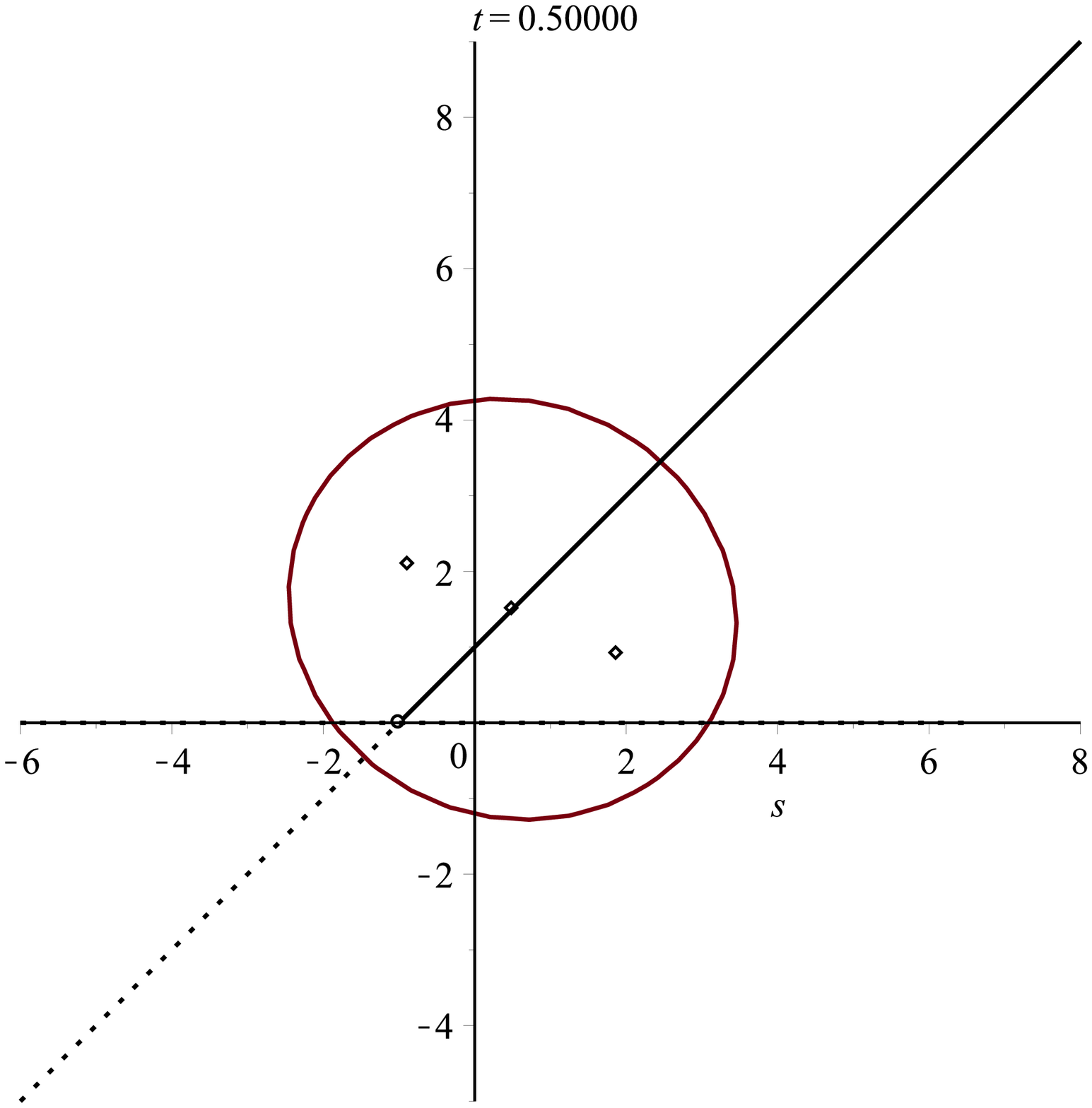}
\includegraphics[scale=0.25]{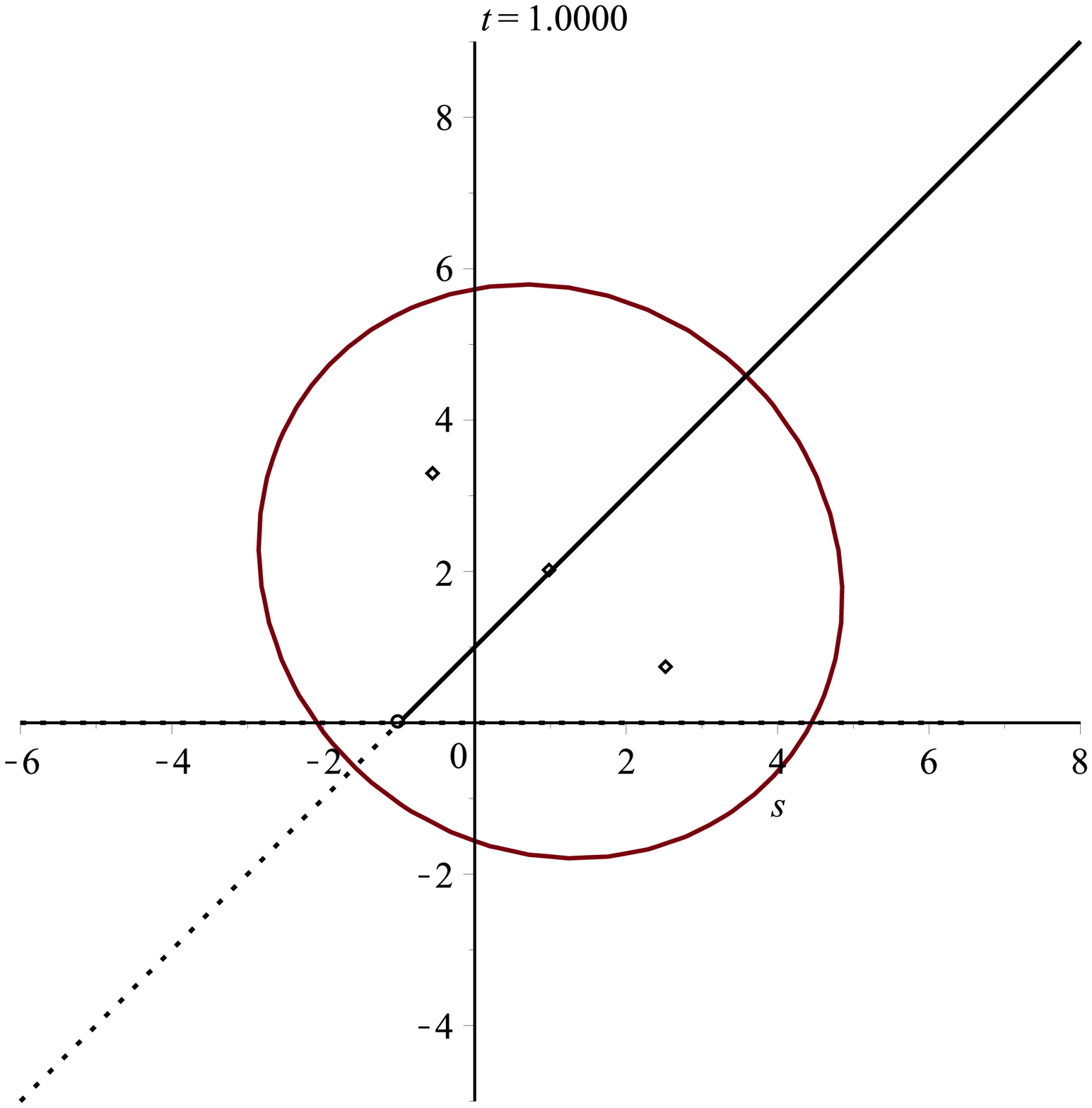}
\caption{Parallel translation along a straight line (from left to right)}
\end{figure}
In case of the straight line the focal set of the translated trifocal ellipse at the parameter $t$ is $-X(t)$, ${\bf 0}$,
$$X(t)= (t+1)\left(\cos \left(\log (t+1) \right)\frac{\partial}{ \partial u^1}_{(t,t+1)}+\sin \left(\log (t+1)\right) \frac{\partial}{ \partial u^2}_{(t,t+1)}\right).$$
Figure 4 shows the parallel translates of the indicatrix along the circle $c(t)=(\cos(t), \sin(t)+2)$. The focal set of the translated trifocal ellipse at the parameter $t$ is $-X(t)$, ${\bf 0}$,
$$X(t)= (\sin(t)+2)\cdot$$
$$\left(\cos \left(\log (\sin(t)+2)\right)\frac{\partial}{ \partial u^1}_{(\cos(t), \sin(t)+2)}+\sin \left(\log (\sin(t)+2)\right) \frac{\partial}{ \partial u^2}_{(\cos(t), \sin(t)+2)}\right).$$
\begin{figure}[h]
\includegraphics[scale=0.25]{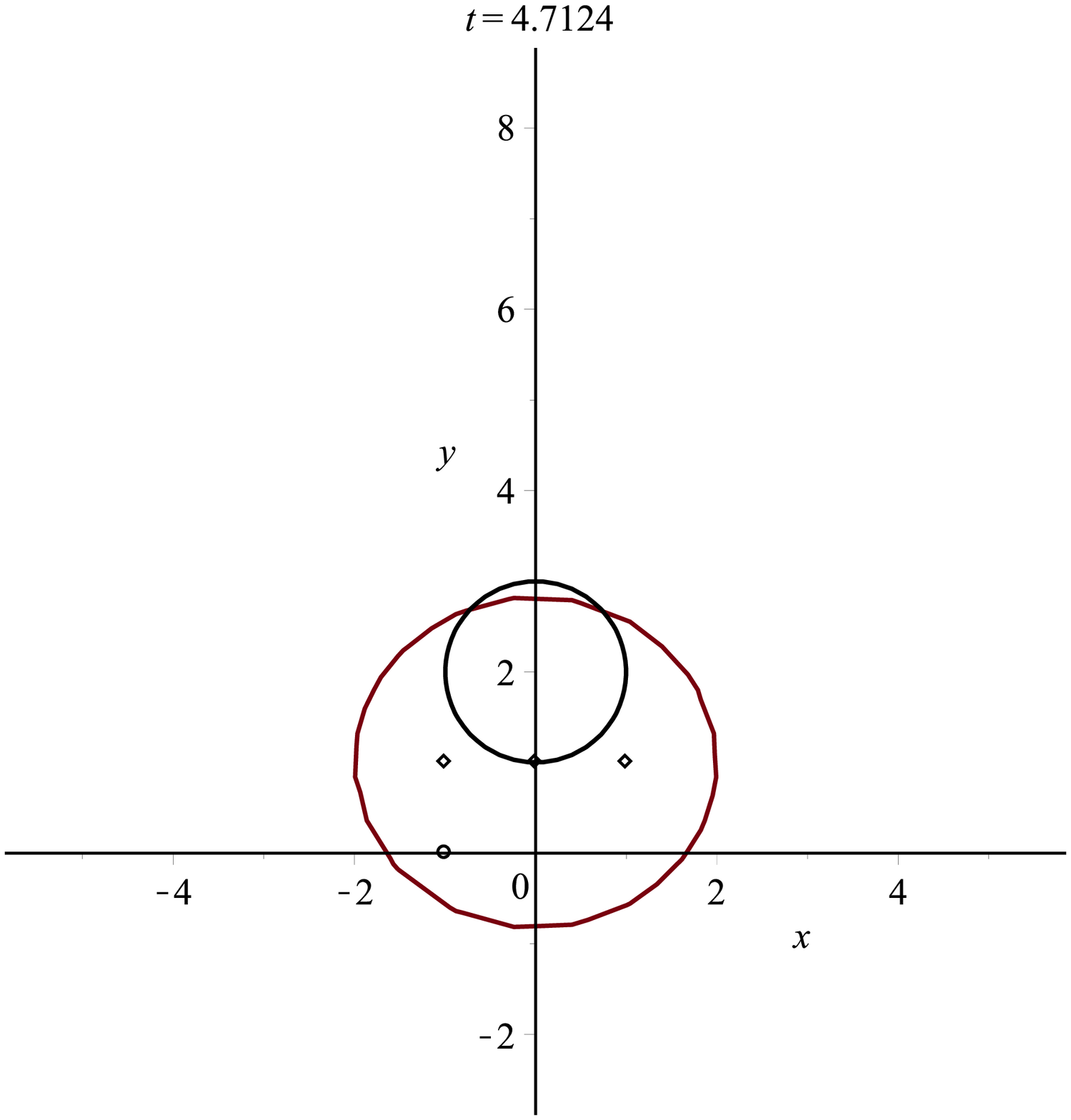}
\includegraphics[scale=0.25]{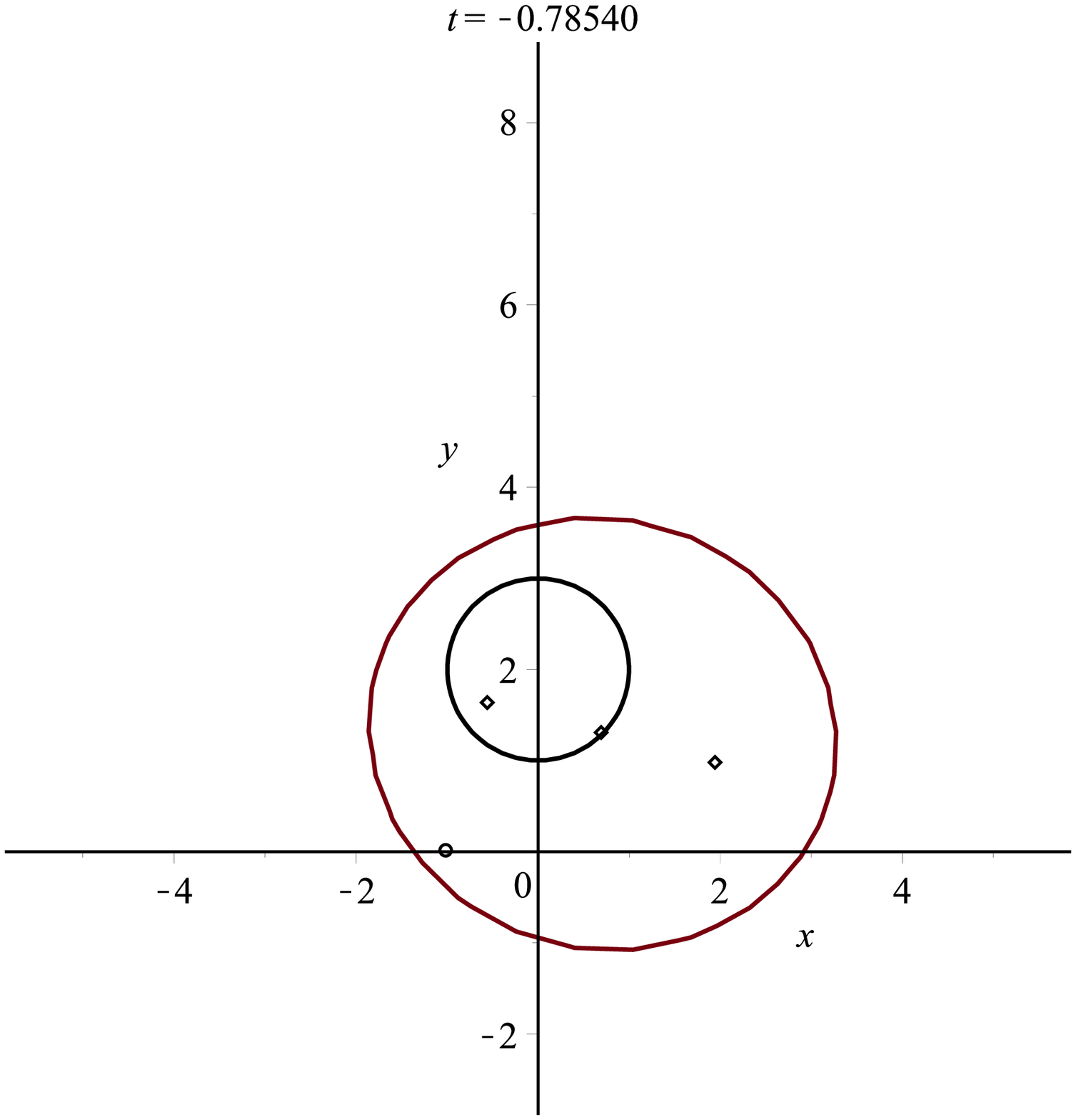}
\includegraphics[scale=0.25]{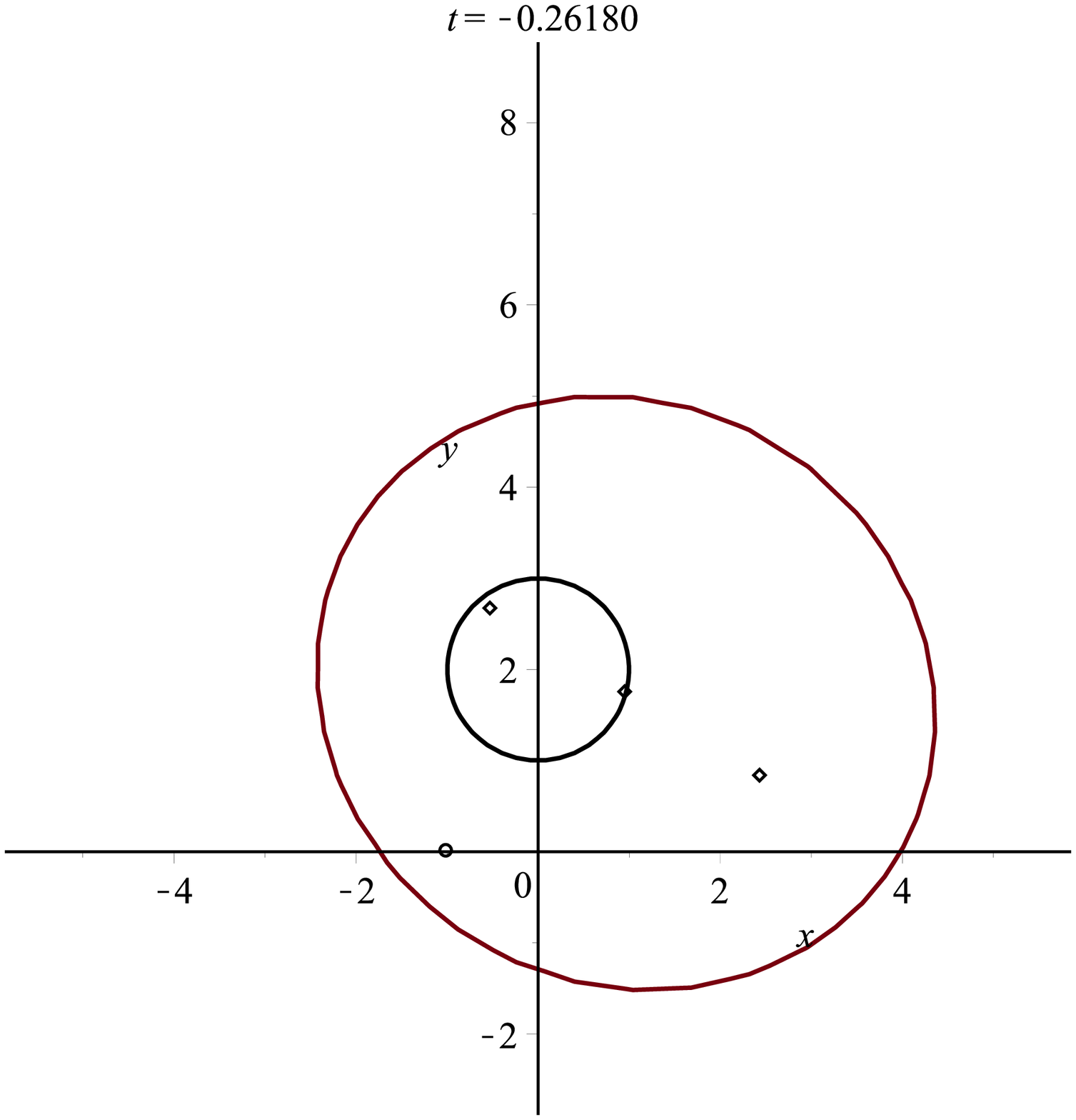}
\end{figure}
\begin{figure}[h]
\includegraphics[scale=0.25]{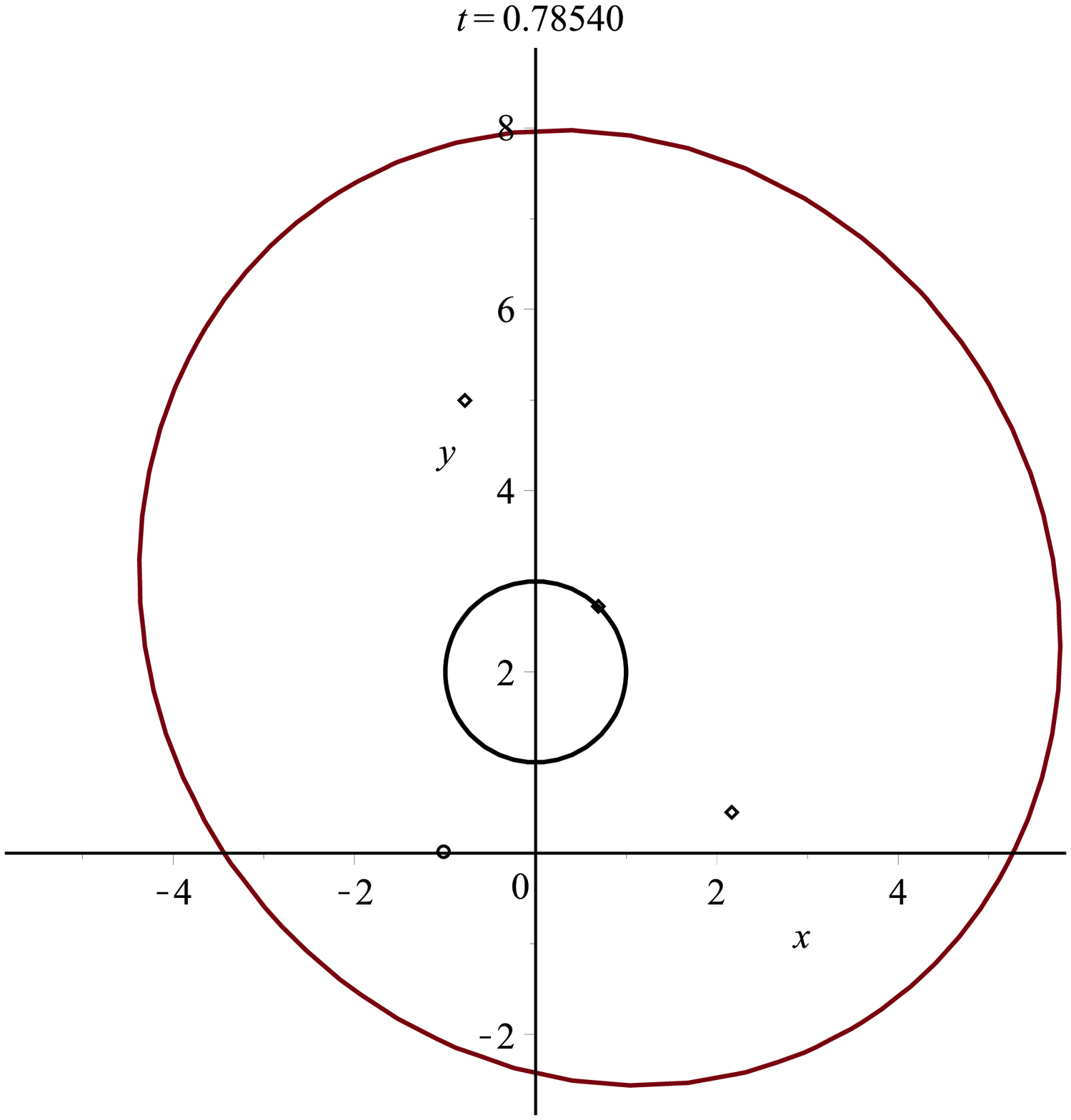}
\includegraphics[scale=0.25]{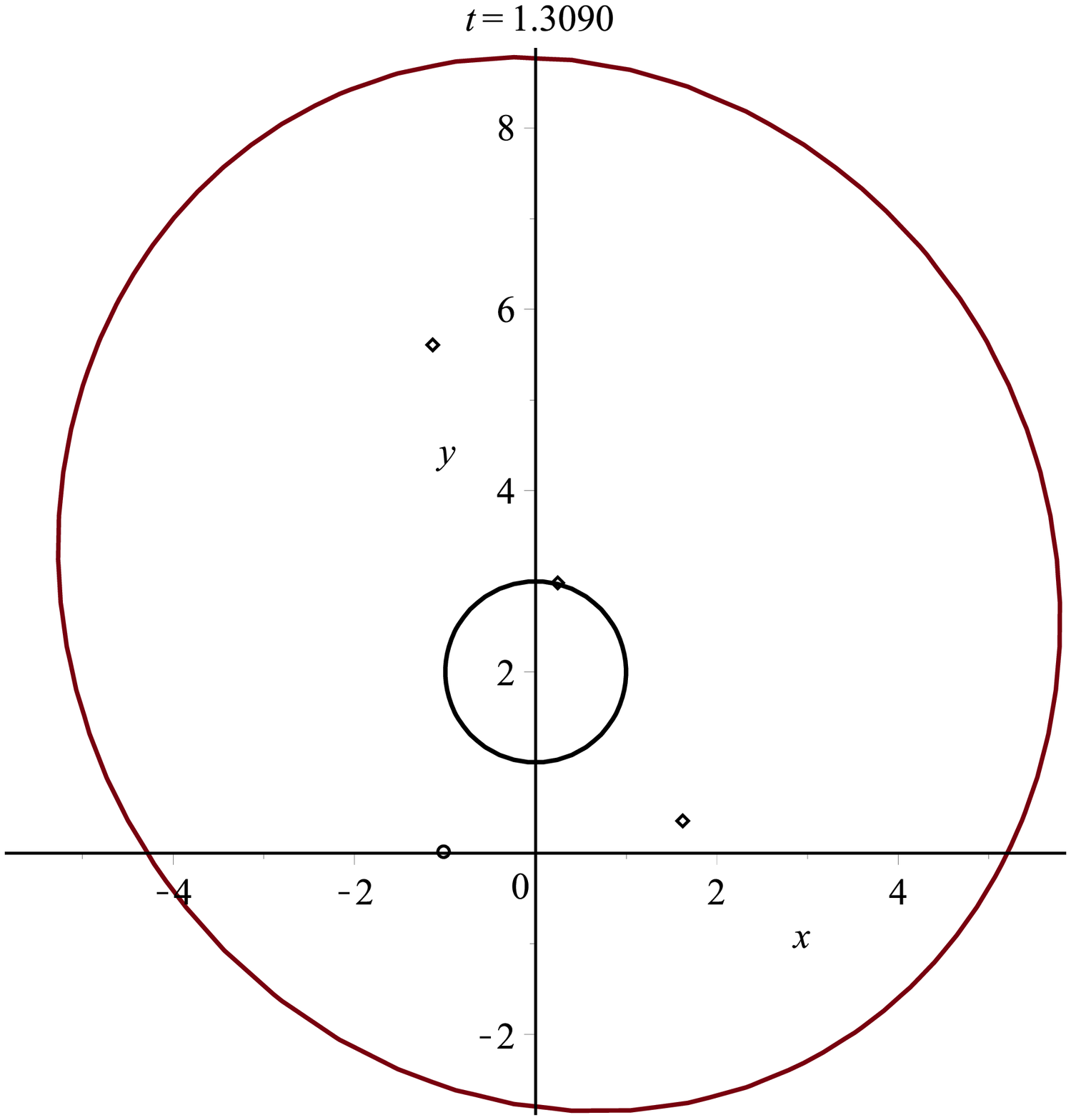}
\includegraphics[scale=0.25]{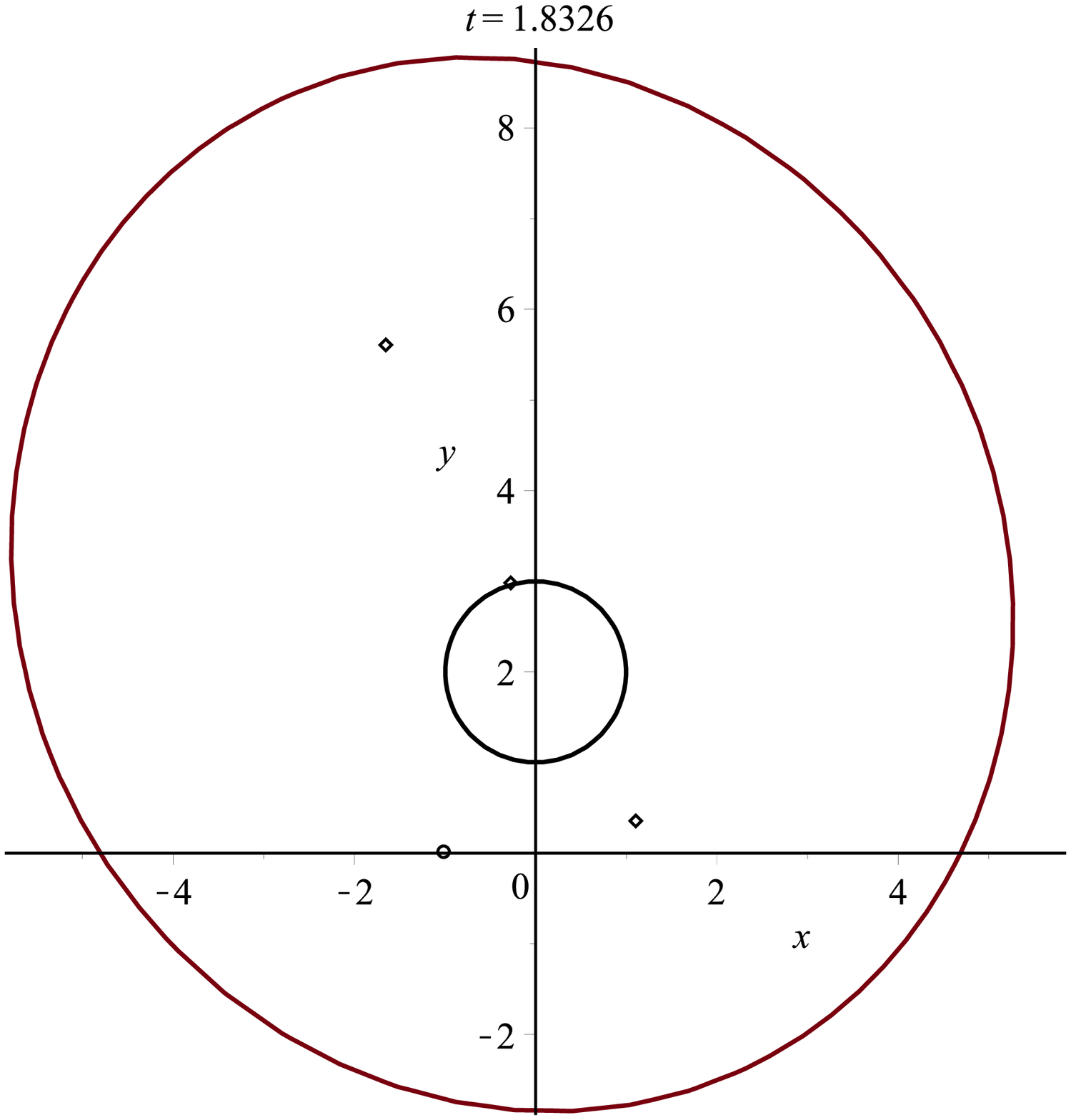}
\end{figure}
\begin{figure}[h]
\includegraphics[scale=0.25]{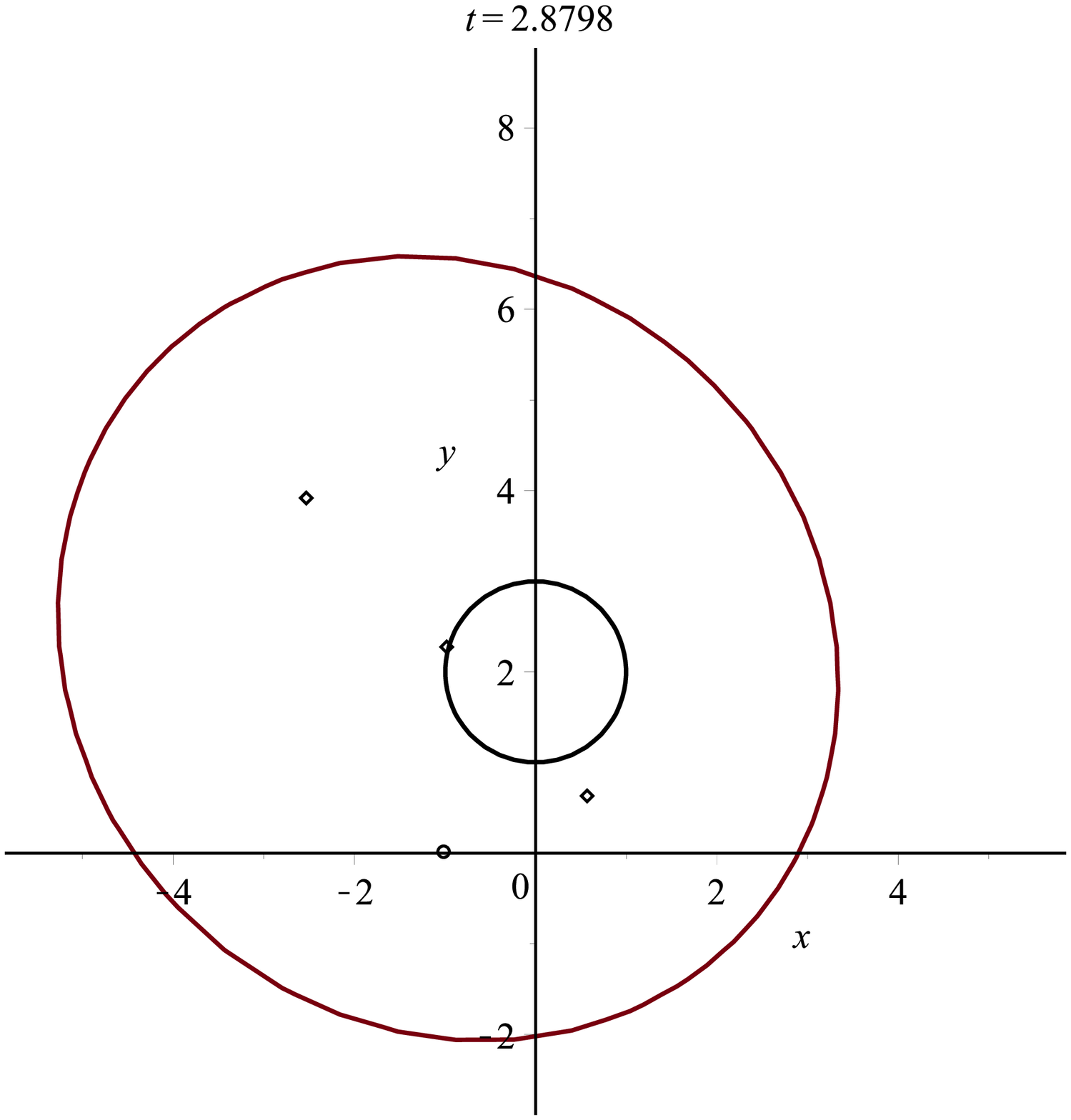}
\includegraphics[scale=0.25]{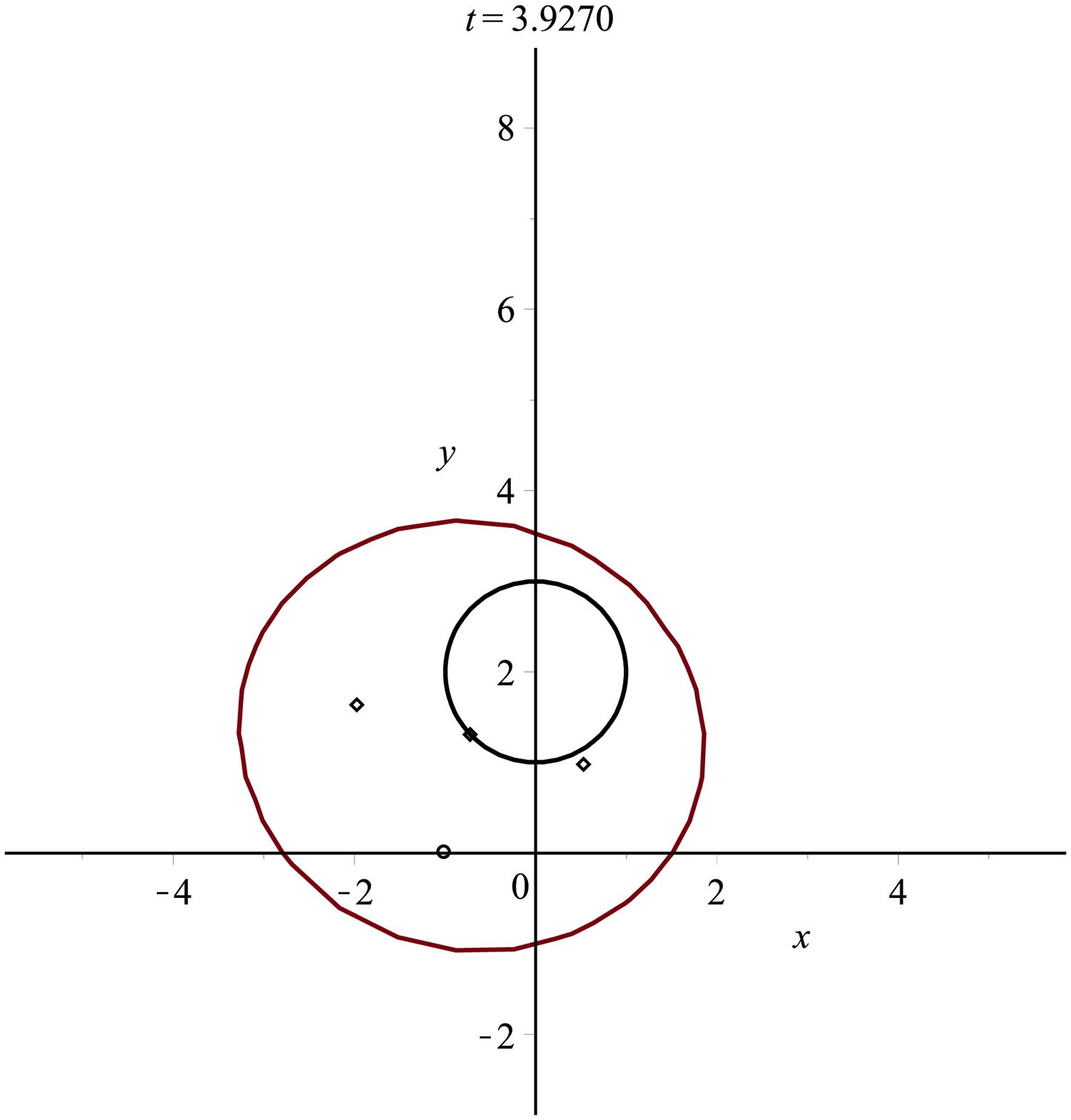}
\includegraphics[scale=0.25]{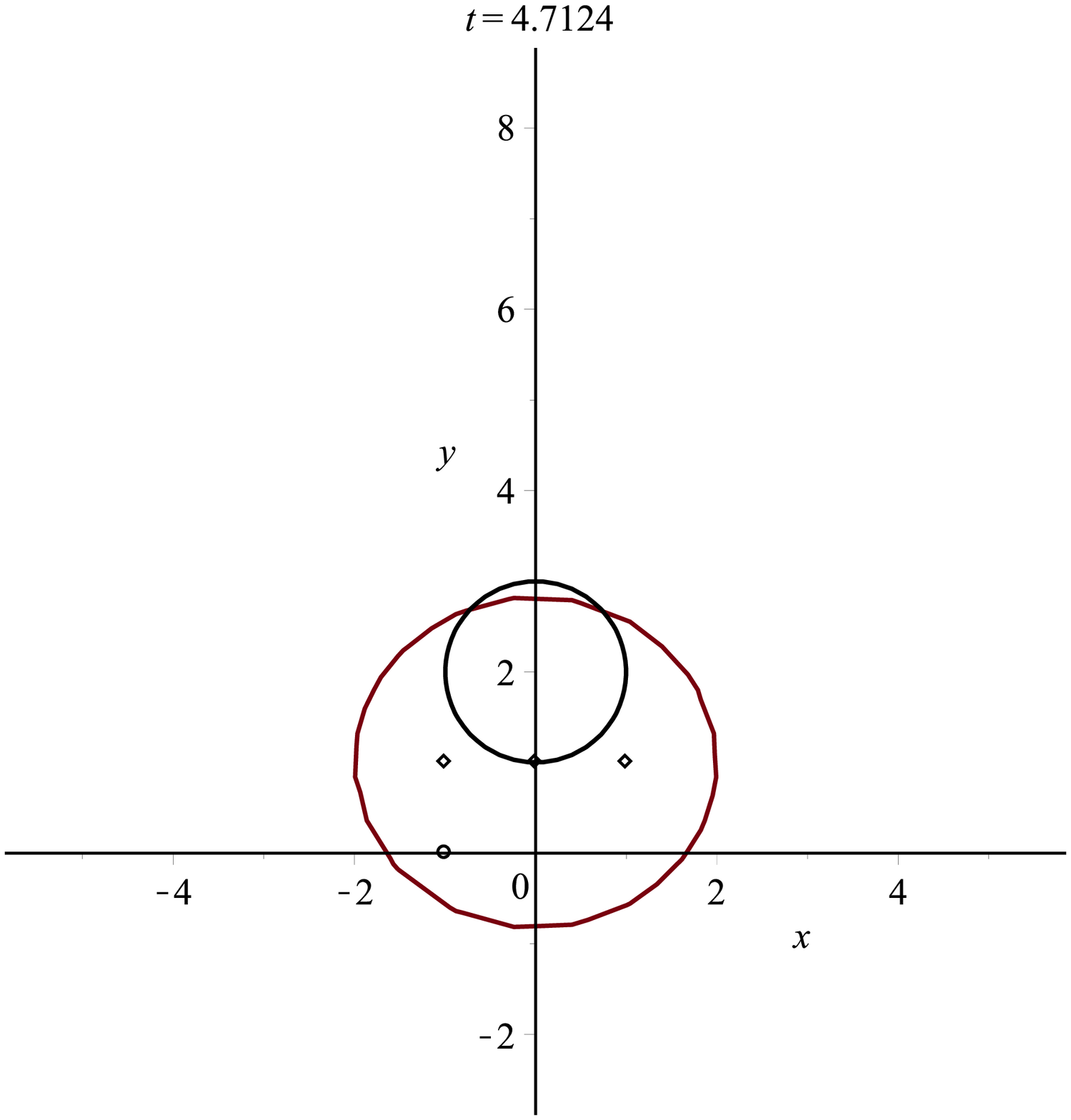}
\caption{Parallel translation along the circle (from left to right)}
\end{figure}

The induced generalized Berwald half plane is not conformally flat. To present examples for conformally flat generalized Berwald manifold it is sufficient and necessary to choose a closed (and, consequently, exact) $1$-form $\rho$.


\begin{thebibliography}{99}

\bibitem{H2}
M. Hashiguchi, \emph{On conformal transformations of Finsler metrics}, J. Math. Kyoto
Univ. 16 (1976), 25-50.

\bibitem{HY}
M. Hashiguchi and Y. Ichijy$\bar{\textrm{o}}$, \emph{On conformal transformations of Wagner
spaces},
Rep. Fac. Sci. Kagoshima Univ. (Math., Phys., Chem.) No. 10 (1977), 19-25.

\bibitem{M2}
M. Matsumoto, \emph {Conformally Berwald and conformally flat Finsler spaces}, Publ.
Math. Debrecen, 58 (1-2)
(2001), 275-285.

\bibitem{V1}
Cs. Vincze, \emph{An intrinsic version of Hashiguchi-Ichijyo's theorems for Wagner manifolds}, SUT J. Math. 35 (2) (1999), 263-270.

\bibitem{V2}
Cs. Vincze, \emph{On Wagner connections and Wagner manifolds}, Acta Math. Hung. 89 (1-2) (2000), 111-133.

\bibitem{V5}
Cs. Vincze, \emph{A new proof of Szab\'{o}'s theorem on the Riemann-metrizability of Berwald manifolds}, J. AMAPN,
21 (2005), 199-204.

\bibitem{V6}
Cs. Vincze, \emph{On a scale function for testing the conformality of Finsler manifolds
to a Berwald manifold}, Journal of Geometry and Physics. 54 (2005), 454-475.

\bibitem{V7}
Cs. Vincze, \emph{On geometric vector fields of Minkowski spaces and their applications}, J. Diff. Geom. and Its Appl. 24 (2006), 1-20.

\bibitem{V9}
Cs. Vincze, \emph{On Berwald and Wagner manifolds}, J. AMAPN,
24 (2008) 169-178.

\bibitem{V10} Cs. Vincze, \emph{Generalized Berwald manifolds with semi-symmetric linear connections}, Publ. Math. Debrecen 83 (4) (2013), pp. 741-755.

\bibitem{V11} Cs. Vincze, \emph{On a special type of generalized Berwald manifolds: semi-symmetric linear connections preserving the Finslerian length of tangent vectors}, European Journal of Mathematics, December 2017, Volume 3, Issue 4, pp 1098 - 1171.

\bibitem{V12} Cs. Vincze, \emph{Lazy orbits: an optimization problem on the sphere}, Journal of Geometry and Physics
Volume 124, January 2018, Pages 180-198. 

\bibitem{Wag1}
V. Wagner, \emph{On generalized Berwald spaces}, CR Dokl. Acad. Sci. USSR (N.S.)
{\bf 39} (1943), pp. 3-5.

\end{thebibliography}
\end{document}